\documentclass[11pt]{article}
\usepackage{amsmath,amsfonts}
\usepackage{graphicx,subfig,multicol}
\usepackage{multirow,rotating,array} 
\usepackage{enumitem} 
\usepackage[top=1.0in, left=1.0in, right=1.0in, bottom=1.0in]{geometry} 
\usepackage[title]{appendix} 
\usepackage{xcolor}

\numberwithin{equation}{section} 

\begin{document}
\centerline{\large \bf Analysis of an Arctic sea ice loss model in the limit of a discontinuous albedo}
\vspace{0.75cm}
\centerline{Kaitlin Hill\footnote{Department of Engineering Sciences and Applied Mathematics, Northwestern University, Evanston, IL 60208, USA (k-hill@u.northwestern.edu)}, Dorian S. Abbot\footnote{Department of the Geophysical Sciences, University of Chicago, Chicago, IL 60637, USA \,(abbot@uchicago.edu)}, and Mary Silber\footnote{Department of Statistics, University of Chicago, Chicago, IL 60637, USA (msilber@uchicago.edu)}}

\vspace{0.5cm}

\noindent{\footnotesize \textbf{Abstract:} As Arctic sea ice extent decreases with increasing greenhouse gases, there is a growing interest in whether there could be a bifurcation associated with its loss, and whether there is significant hysteresis associated with that bifurcation. A challenge in answering this question is that the bifurcation behavior of certain Arctic energy balance models have been shown to be sensitive to how ice-albedo feedback is parameterized. We analyze an Arctic energy balance model in the limit as a smoothing parameter associated with ice-albedo feedback tends to zero, which introduces a discontinuity boundary to the dynamical systems model. Our analysis provides a case study where we use the system in this limit to guide the investigation of bifurcation behavior of the original albedo-smoothed system. In this case study, we demonstrate that certain qualitative bifurcation behaviors of the albedo-smoothed system can have counterparts in the limit with no albedo smoothing. We use this perspective to systematically explore the parameter space of the model. For example, we uncover parameter sets for which the largest transition, with increasing greenhouse gases, is from a perennially ice-covered Arctic to a seasonally ice-free state, an unusual bifurcation scenario that persists even when albedo-smoothing is reintroduced. This analysis provides an alternative perspective on how parameters of the model affect bifurcation behavior. We expect our approach, which exploits the width of repelling sliding intervals for understanding the hysteresis loops, would carry over to other positive feedback systems with a similar natural piecewise-smooth limit, and when the feedback strength is likewise modulated with seasons or other periodic forcing.}
\vspace{0.5cm}


\section{Introduction}
\label{sec:introduction}

A current area of interest in climate science is how perturbations, such as increased greenhouse gases, will affect Arctic sea ice, and what will be the manifestation of this change. One possibility is that the ice, which is currently present year-round, and which cycles between different extremes with the seasons, will cross a ``tipping point'' to  a completely different state. Mathematically, a tipping point may be associated with a bifurcation, after which further parameter changes result in a large and sudden change in system behavior that is not easily reversed (i.e., a hysteresis loop forms). In light of the rapid decrease in summer sea ice extent in the Arctic region during the past decade~\cite{Lenton13}, a potential ``tipping point'' in summer sea ice has been considered in recent studies of conceptual dynamical systems models~\cite{Abbot11,Eisenman09,Eisenman12,Wagner15}. Tipping points have also been extensively investigated in large-scale global climate models, which examine the loss of sea ice from the Arctic region in the coming decades under various greenhouse gas emission scenarios~\cite{Armour11,Holland06,Li13,Tietsche11}.

In this paper we analyze a single-column energy balance model for the Arctic region, originally proposed by Eisenman and Wettlaufer \cite{Eisenman09} (hereafter EW09) as a simple conceptual framework for investigating the role sea ice thermodynamics may play in bifurcations associated with the loss of Arctic sea ice. Their model took the form of a nonautonomous, one-degree-of-freedom ordinary differential equation, which captured the yearly variation in solar forcing, averaged over the Arctic region, as well as seasonal variations in heat transport from lower latitudes into the Arctic. It incorporated the classic ``ice-albedo feedback," a positive feedback mechanism that can lead to coexisting stable states: one state represents a perennially ice-covered Arctic with high albedo (i.e.,~high reflectivity), and the other a perennially ice-free Arctic with low albedo~\cite{North81}. Their model also incorporated the sea ice self-insulation feedback, a thermodynamic feedback mechanism where thicker ice insulates itself better. EW09 used their model to show how sea ice self-insulation, which allows thin ice to grow more rapidly than thick ice, may mitigate, to a degree, positive ice-albedo feedback effects.  Their model shows a smooth transition, under increases in atmospheric greenhouse gases, from a stable perennially ice-covered Arctic to a stable seasonally ice-free state.

The model in EW09 evolves an energy density $E(t)$ that is measured relative to an ocean mixed layer at the freezing point; when $E<0$, $E$ represents the average Arctic sea ice thickness, and when $E>0$ it represents a mean temperature of the Arctic ocean mixed layer. As the energy fluxes change with the state of the system (ice vs.~ocean), this conceptual framework leads naturally to a formulation of the problem as a piecewise-smooth (PWS) continuous dynamical system \cite{Simpson10}: across the ``discontinuity boundary'' at $E=0$, the vector field is Lipschitz continuous but not differentiable. However, in the limit where we take the albedo function to be piecewise-constant, the vector field is generally not continuous across this boundary. Incorporating a piecewise-constant albedo function is a standard approximation exploited in analytic studies of energy balance models first developed to investigate planetary energy balance; in such studies the albedo jumps between high and low values at the latitude of the ice cap boundary~\cite{Budyko69,North75,Sellers69}. In our investigation a boundary in phase space is crossed when the sea ice is lost and then again when it is recovered over the course of the yearly cycle of a seasonally ice-free solution. We focus our analysis on the key role this discontinuity boundary plays in the bifurcation structure of the problem. For a review of PWS systems and bifurcations involving discontinuity boundaries, see \cite{diBernardo}.

EW09 modeled the transition between when the Arctic is ice-covered and ice-free by introducing a smoothing energy scale parameter $\Delta E$ over which the albedo transitions smoothly between a high value (for sea ice) and a low value (for open ocean). The smooth albedo transition is meant to parameterize the spatial variation of ice and capture the albedo dependence on ice thickness \cite{Eisenman09}. However, model results are sensitive to the size of this energy range: transitions between solution states, via pairs of saddle-node bifurcations, are easily ``smoothed out'' by increasing the energy range $\Delta E$ in the albedo transition \cite{Eisenman09}. We take this parameter to zero to formulate our system. An immediate consequence of removing $\Delta E$ from the albedo function is that the resulting system then admits ``sliding regions," which classify it as a Filippov system \cite{Filippov88}. In general, a Filippov system is one for which the flow components normal to the discontinuity boundary at a point may have opposite signs on either side of the boundary. When this occurs a trajectory may ``slide'' along the boundary through so-called sliding regions \cite{diBernardo}. Specifically, if we introduce time as a dynamical variable to render the model of EW09 an autonomous system, then it takes the form
\begin{align}\label{eq:F+F-}
	\frac{dE}{dt} &= \begin{cases}
			G_+(\tau,E), & E>0, \\
			G_-(\tau,E), & E<0, 
		\end{cases} \\
	\frac{d\tau}{dt} &= 1.\label{eq:dtau-dt}
\end{align}
Here $G_+$ and $G_-$ are period-one functions of $\tau$, which is measured in units of years. $E=0$ represents a discontinuity boundary in the $(\tau,E)$-phase plane. Note that the dynamics in the $E=0$ boundary have not been fully specified in the above formulation; this ambiguity in the dynamics is inconsequential for evolving a solution through the discontinuity boundary  {\it provided} that $dE/dt$ has the same sign as $E\to 0^+$ as for $E\to 0^-$. However, for the model of EW09 in the Filippov limit, there are necessarily intervals of $\tau$ for which $G_+(\tau,E)$ and $G_-(\tau, E)$ have opposite signs as $E\to 0^+$ and $E\to 0^-$, respectively. These $\tau$ intervals are associated with ``repelling sliding intervals" in the model, as $dE/dt>0$ for $E\to 0^+$ and $dE/dt<0$ for $E\to 0^-$. Uniqueness of solutions with initial conditions in such sliding intervals is manifestly violated~\cite{diBernardo}. Solutions that enter repelling sliding intervals have the property that they may slide, for an indeterminate time, before either jumping off the discontinuity boundary in either direction or reaching the end of the sliding interval. (Note that the maximum amount of time a solution may slide is equal to the width of the sliding interval, since solutions $E(\tau)$ evolve $\tau$ according to $d\tau/dt=1.$) In this way the lack of uniqueness leads to families of solutions that interact with the sliding intervals. These families of nonunique solutions, when omitted from the analysis, create ``gaps'' in the bifurcation diagram where solution branches become discontinuous.

Our analysis contributes a new case study of how the Filippov limit of a dynamical system, where the vector field is discontinuous and which here is easier to analyze in closed form, can be utilized to determine properties of the original continuous system. The qualitative comparison between a smooth system and its nonsmooth counterpart has been of significant recent interest in the context of applications \cite{Harris15,Machina13} as well as in general mathematical investigations that consider either ``smoothing out'' a PWS system \cite{Jeffrey14,Kristiansen15,Teixeira12} or transforming a smooth system into a PWS one \cite{Desroches11}. We explore how the relative widths of the sliding intervals affect the bifurcation structure, which provides insight into how model parameters impact the transition from an ice-covered to an ice-free Arctic. For example, we show that a small hysteresis loop, involving seasonally ice-free states and associated with the loss of perennial ice, can be moved along the seasonal solution branch by varying the relative widths of the sliding intervals. In particular, it can be moved so that it is associated with the loss of the perennially \textit{ice-free} Arctic state rather than the perennially \textit{ice-covered} Arctic state. Additionally, we use insight gained from this perspective to determine regions of parameter space where the transition from an ice-covered to a seasonally ice-free Arctic is large.

This paper is organized as follows. Section \ref{sec:model} describes the original formulation of the Eisenman and Wettlaufer model \cite{Eisenman09}, as well as variations on it \cite{Abbot11,Eisenman12}, and introduces our general Filippov formulation of the model. In Section \ref{sec:solutions} we derive, in turn, the existence and stability conditions for the three types of periodic solutions: perennially ice-free, perennially ice-covered, and seasonally ice-free. Section \ref{sec:case-study} presents a case study using the model given in \cite{Eisenman12}. In this case study we compare numerical bifurcation diagrams for the albedo-smoothed system with analytical bifurcation diagrams for the corresponding Filippov system, and we explore the various parameters and how they impact the bifurcation structure in an alternative manner to previous parameter studies of the model. A summary and discussion of our results concludes the paper in Section \ref{sec:discussion}.

\section{Model Description}
\label{sec:model}

The model developed in EW09 is an energy balance model, where the energy per unit surface area, $E$, is measured relative to  an Arctic ocean mixed layer at the freezing point.  An increase in the energy density relative to this baseline leads to warming of the mixed layer to a temperature $T\propto E$ above the freezing point, while a decrease in the energy density leads to sea ice growth to an ice thickness $h_i\propto -E$. Specifically, the state variable $E$ in EW09 has the physical interpretation that
\begin{equation}\label{eq:Edef}
	E=\begin{cases}
	-L_ih_i, & E<0,\\
	C_sT, & E\ge 0,\ \end{cases}
\end{equation}
where the proportionality constants, $L_i$ and $C_s$, correspond to the latent heat of fusion for sea ice and the ocean heat capacity per unit surface area, respectively.   

The evolution equation for $E(t)$ in EW09 is derived by considering heat flux balances at the surface, and it takes the form \cite{Eisenman09}
\begin{equation}\label{eq:EW09}
	\frac{dE}{dt}=\left(1-\alpha(E)\right)F_S(t) - \left(F_0(t)+F_T(t)T(t,E)\right) + F_B - \nu_0 {\cal H}(-E)E.
\end{equation}
The energy source term $(1-\alpha(E))F_S(t)$ represents the absorbed incoming radiation, where $F_S(t)$ is the solar radiation flux, averaged over the Arctic region, and $(1-\alpha(E))$ is the fraction of the incoming radiation absorbed ($\alpha(E)$ is the albedo). The energy loss term $(F_0(t)+F_T(t)T(t,E))$ represents outgoing radiation, which increases with the surface temperature $T(t,E);$ it has been linearized about the freezing point $T=0$. The flux terms $F_S(t)$, $F_0(t)$, and $F_T(t)$ are all periodic with period 1 year. Here $F_B$ represents a small, constant basal heat flux into the Arctic mixed layer/sea ice from the deep ocean. The term $\nu_0 {\cal H}(-E)E$ captures sea ice transport out of the Arctic, where ${\cal H}(x)$ is the Heaviside function; this term is only present if there is sea ice. In EW09, the albedo function smoothly transitions from its ice value ($\alpha_i$) to its ocean mixed layer value ($\alpha_{ml}$) and is given by
\begin{equation}\label{eq:AEW09}
	\alpha(E) = \frac{\alpha_{ml}+\alpha_i}{2} + \frac{\alpha_{ml}-\alpha_i}{2}\tanh\left(\frac{E}{\Delta E}\right),
\end{equation}
where $\Delta E$ is the previously mentioned energy scale for smoothing.

Sea ice thermodynamics enters the model in EW09 by the surface temperature function $T(t,E),$ which depends on the state of the ice/mixed layer:
\begin{equation*}
	T(t,E) = \begin{cases}
		E/C_s, &\qquad E\ge0, \\
		0,     &\qquad E<0, \quad (1-\alpha_i)F_S(t)>F_0(t), \\
		\left(\frac{(1-\alpha_i)F_S(t)-F_0(t)}{F_T(t)}\right)\left(\frac{E}{E-k_i L_i/F_T(t)}\right), &\qquad E<0, \quad (1-\alpha_i)F_S(t)\leq F_0(t),	
	\end{cases}
\end{equation*}
where $k_i$ is the thermal conductivity of sea ice. Note that when $E\geq0$, no sea ice is present, so $T$ is determined by (\ref{eq:Edef}). When $E<0$, sea ice is present, and $T$ depends on whether or not the ice is melting at the surface. When the outgoing radiation exceeds the incoming radiation ($(1-\alpha_i)F_S(t) \leq F_0(t)$), surface temperature depends on the thickness of the sea ice, $h_i \propto -E$, and it is determined, in a large Stefan number approximation, from surface flux balance, $k_i T/h_i = (1-\alpha_i)F_S(t)-F_0(t)-F_T(t)T.$ (See \cite{Eisenman09} for more details.) On the other hand, when $E<0$ and incoming exceeds outgoing radiation, the ice is melting, and in a large Stefan number approximation the surface temperature immediately warms to the melting point, defined here to be $T=0$.
 
A variation of the model in EW09, due to Abbot, Silber, and Pierrehumbert~\cite{Abbot11}, replaces the seasonally varying $F_0(t)$ with a state-dependent model and incorporates possible effects due to cloud cover into the albedo function. In another variation, Eisenman \cite{Eisenman12} simplified the EW09 model by introducing a sinusoidal approximation for each of the seasonally varying terms $F_S(t)$ and $F_0(t)$, while neglecting any seasonal variations of $F_T(t)$ in (\ref{eq:EW09}). This is the version of the model that we use in our case study in Section \ref{sec:case-study}.

Our analysis is based on these Arctic sea ice models in the following form:
\begin{align}\label{eq:dEdt}
	\begin{aligned}
		\frac{dE}{dt}&=F(\tau,E)-BT(\tau,E),\\
		\frac{d \tau}{dt}&=1,
	\end{aligned}
\end{align}
where after suitable nondimensionalization of $E$ \cite{Eisenman12}, the surface temperature is given by
\begin{equation}\label{eq:TE12}
	T(\tau,E)=\begin{cases}
		E, &\qquad E\ge 0, \\
		0, &\qquad E<0, \quad F>0, \\
		\frac{F}{B}\left(\frac{E}{E - \zeta}\right), &\qquad E<0, \quad F \le 0,	
	\end{cases}
\end{equation}
where $\zeta>0$ is a dimensionless thermodynamic parameter associated with sea ice. This is the form of (\ref{eq:EW09}) from EW09 in the case that basal heat transport and sea ice transport out of the Arctic are neglected, i.e.,~$F_B=\nu=0.$ Moreover, we have neglected any temporal variation of $F_T(t)$ in (\ref{eq:EW09}) by replacing it with the constant $B$. (Similar approximations were made in some numerical simulation results in \cite{Eisenman12}.) Note that although $T(\tau,E)$ is piecewise-defined, it is Lipschitz continuous, so it does not contribute to the discontinuity in the vector field at $E=0.$

In the limit of no albedo-smoothing, $F(\tau,E)$ in (\ref{eq:dEdt}) is given by
\begin{equation}\label{eq:F}
	F(\tau,E) = \left\{ 
	\begin{array}{l   l}
		 (1 + \Delta \alpha){\cal F}_s(\tau)-{\cal F}_{l,+}(\tau) \equiv F_+(\tau),  &  E>0, \\
		 (1 - \Delta \alpha){\cal F}_s(\tau)-{\cal F}_{l,-}(\tau) \equiv F_-(\tau),  &  E<0, \\
	\end{array} \right.
\end{equation}
where $F_\pm(\tau)$ are (nondimensionalized) period-one functions of $\tau$. The albedo function here is the limit as $\Delta E\rightarrow 0$ of the albedo function after nondimensionalization in \cite{Eisenman12}; the term
\begin{equation}\label{eq:Da}
	\Delta \alpha\equiv \frac{\alpha_i - \alpha_{ml}}{2\left( 1- (\alpha_i + \alpha_{ml})/2 \right)}\in [0,1]
\end{equation}
represents the difference between the ice and ocean mixed layer albedos, scaled as part of the nondimensionalization process by twice the mean coalbedo (coalbedo = 1 -- albedo), which is approximately 1. The fluxes ${\cal F}_s(\tau)$ and ${\cal F}_{l,\pm}(\tau)$ characterize incoming shortwave and outgoing longwave radiation terms, respectively, akin to $F_S(t)$ and $F_0(t)$ in (\ref{eq:EW09}). The subscript $\pm$ on the outgoing longwave radiation term indicates that it may depend on whether there is ice in the Arctic or not (e.g.,~as in the variations of the model of EW09 investigated by Abbot, Silber, and Pierrehumbert~\cite{Abbot11}). In our case study in Section \ref{sec:case-study} we use
\begin{equation}
	F(\tau,E) = 
	\begin{cases}
		F_+(\tau) \equiv (1 + \Delta\alpha)(1-S_a \cos(2 \pi \tau)) - (L_m + L_a \cos(2 \pi (\tau-\phi))), & E>0, \\
		F_-(\tau) \equiv (1 - \Delta\alpha)(1-S_a \cos(2 \pi \tau)) - (L_m + L_a \cos(2 \pi (\tau-\phi))), & E<0,
	\end{cases}
\end{equation}
which was generated by taking the limit as $\Delta E \rightarrow 0$ of $F(\tau,E)$ from \cite{Eisenman12},
\begin{equation}\label{eq:FE12smooth}
	F(\tau,E) = (1 + \Delta\alpha\tanh(E/\Delta E))(1-S_a \cos(2\pi\tau)) - (L_m + L_a \cos(2\pi(\tau-\phi))).
\end{equation}
This represents a sinusoidal approximation to the periodic functions $F_S(t)$ and $F_0(t)$ from EW09 \cite{Eisenman12}. We note that $L_m,$ the annual mean outgoing longwave radiation when $T=0$, decreases as greenhouse gas concentrations increase \cite{Eisenman12} and is the parameter we vary in our bifurcation studies.

\begin{figure}[!t]
	\captionsetup[subfloat]{farskip=-1.25ex,captionskip=-0.5ex} 
	\centering
	\begin{tabular}{cc}
    \multirow{2}*[2.075in]{\subfloat[]{\includegraphics[scale=1]{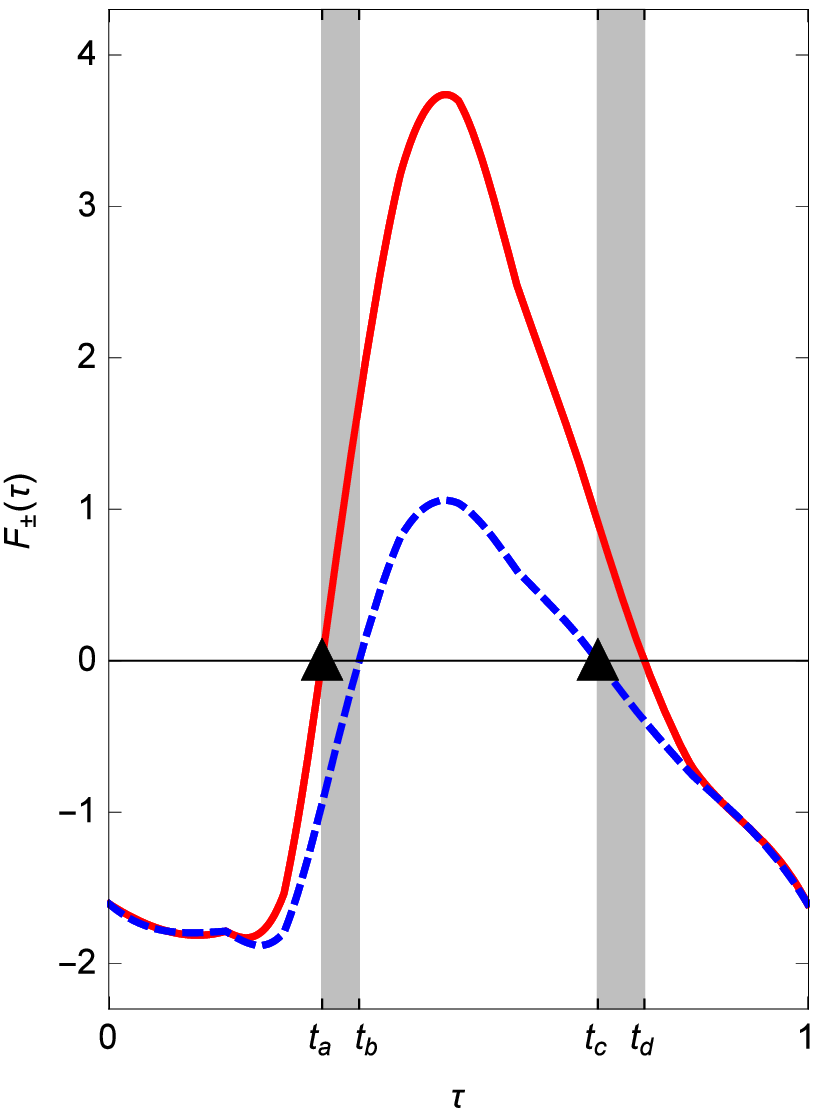}}} 
		  & \subfloat[]{\includegraphics[scale=1]{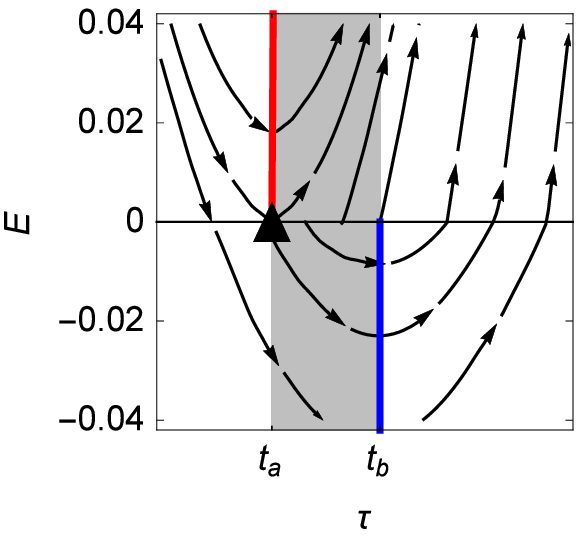}} \\
	    & \subfloat[]{\includegraphics[scale=1]{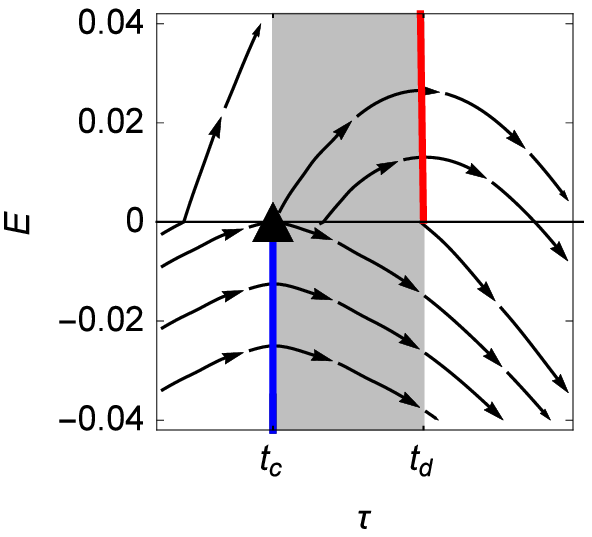}} \\
  \end{tabular}
	\caption{(a) Sample plot of $F_+(\tau)$ (solid red curve), which applies for $E > 0$, and $F_-(\tau)$ (dashed blue curve) for $E<0$, from (\ref{eq:F}) over one period of the seasonal forcing in EW09~\cite{Eisenman09}, after nondimensionalization. Repelling sliding intervals are shaded; triangles indicate the point where a trajectory may enter a sliding interval. For perennially ice-free solutions, $t_a$ is when the minimum of $E(t)$ occurs as $E\rightarrow 0^+$. For perennially ice-covered solutions, ice ablation initiates and terminates at $t_b$ and $t_c$, respectively. For seasonally ice-free solutions, we show in Section~\ref{sec:solutions} that $t_a$ is the last possible freezing time and $t_c$ is the last possible melting time. (b),(c) Blow-up of the direction field for the nondimensionalized EW09 model near the sliding intervals; nullclines are given in red for $E>0$ and blue for $E<0.$ Shaded regions represent the sliding intervals $S_1\equiv[t_a,t_b],$ associated with ``spring,'' and $S_2\equiv[t_c,t_d],$ associated with ``autumn.''}
	\label{fig:fplusfminus}
\end{figure}

The general system we consider can be summarized by combining (\ref{eq:dEdt})--(\ref{eq:F}):
\begin{equation}\label{eq:system-summary}
	\begin{aligned}
		\frac{dE}{dt}&=\begin{cases}
				F_+(\tau) - B E(\tau), &\qquad E > 0, \\
				F_-(\tau), &\qquad E<0, \quad F_->0,\\
				\zeta F_-(\tau)/\left(\zeta - E(\tau)\right), &\qquad E<0, \quad F_- \le 0,
			\end{cases} \\
		\frac{d \tau}{dt}&=1.
	\end{aligned}
\end{equation}
We will refer to this system as the ``Filippov system,'' and we will refer to the original system with albedo smoothing, i.e.,~$\Delta E \neq 0$ in (\ref{eq:FE12smooth}), as the ``albedo-smoothed system.'' There are three possible discontinuity boundaries in this system: one at $E=0$ that is a result of the albedo jump, and two corresponding to roots of $F_-(\tau),$ which are associated with changes in the thermodynamics between growing and melting sea ice. Across each $F_-=0$ boundary the slope of solution trajectories jumps for $E<0$ solutions. For our purposes, though, these boundaries associated with $F_-=0$ are benign, in the sense that trajectories cross them with transverse speed given by $d\tau/dt=1.$ In contrast, the discontinuity boundary $E=0$ associated with the albedo jump may introduce sliding regions; our analysis will focus on how solutions interact with these. Solution behavior on either side of the discontinuity boundary at $E=0$ is characterized through $F_\pm(\tau)$, given by (\ref{eq:F}). The time intervals for which $F_+(\tau)\ge 0$ and $F_-(\tau)\le 0$ correspond to repelling sliding intervals for the Filippov dynamical system (\ref{eq:system-summary}). These time intervals, one in spring and one in autumn, are indicated by shading in Figure \ref{fig:fplusfminus}. Note that the thermodynamic parameter $\zeta$ and the parameter $B$ from (\ref{eq:system-summary}) do not enter into determining the intervals of sliding. The discontinuity boundaries meet at values where $E=0$ and $F_-=0$ (e.g., the points $(t_b,0)$ and $(t_c,0)$ in Figure \ref{fig:fplusfminus}). Since the $F_-=0$ boundaries are ``crossing regions'' they do not interfere with the crossing or sliding behavior on the $E=0$ boundary at these points.

Note that the only way for a solution trajectory to enter a repelling sliding interval is if $E=0$ at the start of a sliding interval (e.g.,~if $E\to 0^+$ as $\tau \to t_a$ or if $E\to 0^-$ as $\tau \to t_c$ in Figure \ref{fig:fplusfminus}). A technique called ``Filippov's convex method'' can be used to define the behavior of the solution in the boundary when it is sliding \cite{Filippov88}. This method can in general reveal diverse phenomena such as equilibria of the sliding flow \cite{diBernardo}, but in our case it produces trivial sliding motion $\langle dE/dt,d\tau/dt\rangle=\langle 0,1\rangle$ due to the system dynamics $d\tau /dt=1.$ 

We do not consider the possibility of ``attracting sliding intervals" for this model since the primary source of discontinuity in the vector field is due to the albedo jump $\Delta \alpha>0$ at $E=0$, i.e.,~through the contribution $(1\pm\Delta \alpha){\cal F}_s(\tau)\geq 0$ to (\ref{eq:F}). This contribution is enough to ensure that $F_+(\tau)\geq F_-(\tau)$ over most of the yearly cycle of forcing. An attracting sliding interval would require that $dE/dt < 0$ for $E\to 0^+$ and $dE/dt > 0$ for $E\to 0^-$, which would require $F_+(\tau)<F_-(\tau)$ over some $\tau$-interval. Attracting sliding intervals do occur, for some parameter sets, in the Filippov limit of the simplified model investigated in \cite{Eisenman12}, but they are a nonphysical consequence of the sinusoidal approximations used for the periodic forcing functions. For this reason we focus our bifurcation analysis in Section \ref{sec:case-study} on parameter sets for which attracting sliding does not occur.

\section{Nonsliding Periodic Solutions}
\label{sec:solutions}

In this section we determine the existence and stability properties of the three distinct types of periodic solutions of the Filippov system (\ref{eq:system-summary}):
\begin{enumerate}[nolistsep]
	\item[A.] Perennially ice-free ($E(\tau)> 0$ for all $\tau$).
	\item[B.] Perennially ice-covered ($E(\tau)<0$ for all $\tau$).
	\item[C.] Seasonally ice-free ($E(\tau)$ takes on both positive and negative values).
\end{enumerate}
By virtue of their definitions the perennially ice-free ($E>0$) and perennially ice-covered ($E<0$) solutions avoid the discontinuity boundary $E=0$. Nonetheless, the parameter space existence boundary for these solutions is determined by the condition that the minimum (maximum) value of $E$ limits to $E=0$. This limiting minimum value of $E=0$ occurs for perennially ice-free solutions at $\tau=t_a$, while the maximum value of $E=0$ occurs for perennially ice-covered solutions at $\tau=t_c$; see Figure~\ref{fig:fplusfminus}(b),(c). In contrast the seasonally ice-free solutions necessarily cross the discontinuity boundary $E=0$ twice during their yearly cycle. In Section~\ref{sec:seasicefree} we show that in addition to saddle-node bifurcations along the seasonal solution branch, these states may be destroyed, as parameters vary, when they collide with the sliding intervals at either $(\tau,E)=(t_a,0)$ or $(\tau,E)=(t_c,0)$.

\subsection{Perennially Ice-Free Periodic Solutions}
\label{sec:pericefree}

In this section we derive conditions for the existence, stability, and uniqueness of perennially ice-free periodic solutions that do not enter the sliding interval. Specifically, we show that these solutions are always unique and stable when they exist, and that this solution branch arises via a ``grazing-sliding bifurcation'' \cite{diBernardo} with the discontinuity boundary $E=0$ at $\tau=t_a$. A grazing-sliding bifurcation occurs when a trajectory that resides entirely on one side of the discontinuity boundary is continuously perturbed into a trajectory that grazes the boundary at a point and, with further perturbation, may begin to slide \cite{diBernardo}.

Perennially ice-free solutions satisfy $E > 0$ for all $\tau$, so that based on (\ref{eq:system-summary}), $E$ evolves according to  
\begin{equation}\label{eq:icefree}
	\frac{dE}{d\tau}=F_+(\tau)-BE.
\end{equation}
(For Section \ref{sec:solutions} we use the notation $\frac{dE}{d\tau}$ for simplicity, since $\frac{dt}{d\tau}=1$.) From this we can derive a period-one Poincar\'e return map, which takes an initial condition $E_0\equiv E(t_0)$ at time $t_0\in[0,1)$ to its value one period later:
\begin{equation}
	E(t_0+1)=e^{-B}[E(t_0)+e^{-Bt_0}I_+(t_0,t_0+1)],
\end{equation}
where
\begin{equation}\label{eq:I+}
	I_+(t_0,\tau)\equiv\int_{t_0}^{\tau} e^{Bt}F_+(t)dt.
\end{equation}
(At this stage we have not chosen a preferred phase $t_0$ for our map; it is chosen below.) Periodic solutions satisfy $E(t_0+1)=E(t_0)\equiv E_0^*$, which determines the unique initial condition for the periodic state:
\begin{equation}\label{eq:icefreeic}
	E_0^*=\left(e^B-1\right)^{-1}e^{-Bt_0}I_+(t_0,t_0+1).
\end{equation}
The stability of this periodic state is determined by the eigenvalue of the linearized Poincar\'e map, which is the Floquet multiplier 
\begin{equation}
	\mu=dE(t_0+1)/dE(t_0)=e^{-B}.
\end{equation}
Since $B$ is a positive constant, $|\mu|<1$ and the periodic solution is always stable when it exists.

Finally, we determine the existence condition for perennially ice-free solutions, which is that $E(\tau )>0$ for all values of $\tau\in[t_0,t_0+1)$. To do this we consider the time at which $E$ achieves its minimum, $t_{min},$ and let $E(t_{min})\to 0$. The conditions for a minimum, that $\frac{dE}{d\tau}=0$ and $\frac{d^2E}{d\tau^2}>0$, become simply $F_+(t_{min})=0$ and $\frac{dF_+}{d\tau}>0$ for $E(t_{min})\rightarrow 0,$ so $t_{min}\rightarrow t_a$. (See Figure~\ref{fig:fplusfminus}(a) for a sample plot of $F_+(\tau)$.) Hence we may set $t_0=t_{min}=t_a$ in (\ref{eq:icefreeic}) so that perennially ice-free solutions exist provided $E_0^*>0$, i.e.,~provided
\begin{equation}\label{eq:ice-free-exist}
	I_+(t_a,t_a+1)>0,
\end{equation}
where $t_a$ is determined from
\begin{equation}
	F_+(t_a)=0,\quad \frac{dF_+}{d\tau}(t_a)>0.
\end{equation}

In this analysis of the perennially ice-free periodic solutions we have avoided the sliding intervals of the Filippov system by insisting that $E(\tau)>0$ for all $\tau$. However, the existence boundary for these solutions in parameter space corresponds to conditions where $E(t_a)\rightarrow 0$, which is also the unique entry point to $S_1$, the springtime sliding interval in Figure \ref{fig:fplusfminus}(b). Thus the bifurcation that leads to the creation of the perennially ice-free solutions corresponds to a grazing-sliding bifurcation \cite{diBernardo}. In relation to the classification of grazing-sliding bifurcations in \cite{Kuznetsov03}, this bifurcation is similar to the generic case $TC_2,$ where a sliding periodic trajectory and a nonsliding periodic trajectory collide and annihilate each other.

\subsection{Perennially Ice-Covered Periodic Solutions}
\label{sec:pericecov}

For these periodic states, $E(\tau)<0$ for all $\tau$. We introduce a time $t_b\in[0,1)$ when surface ablation initiates so that $T=0$ in~(\ref{eq:dEdt}), and another time $t_c\in(t_b,t_b+1)$ when surface temperatures next drop below the freezing point. Our focus is on periodic solutions, and we consider those over the time interval $\tau\in[t_b,t_b+1]$. There are thus two phases associated with this solution: a summer phase, $\tau\in[t_b,t_c],$ when ice is ablating from the surface and the surface temperature is fixed at $T=0,$ and a winter phase when $T$ drops below the freezing point, which occurs for $\tau\in(t_c,t_b+1).$ Specifically, $t_b$ and $t_c$ are defined by the conditions
\begin{align}\label{eq:tbtcdef}
	\begin{aligned}
		F_-(t_b)&=0,\quad \frac{dF_-}{d\tau}(t_b)>0,\\
		F_-(t_c)&=0,\quad \frac{dF_-}{d\tau}(t_c)<0,
	\end{aligned}
\end{align}
as shown in Figure \ref{fig:fplusfminus}(a). Based on these two phases, (\ref{eq:system-summary}) simplifies to
\begin{equation}
	\frac{dE}{d\tau} = \left\{ 
		\begin{array}{ll}
			F_-(\tau), &  \tau\in[t_b,t_c], \\
			\frac{\zeta F_-(\tau)}{\zeta - E},& \tau\in(t_c,t_b+1), \\
		\end{array} \right. 
\end{equation}
Note that the times $t_b$ and $t_c$ separate the phases when the perennial ice is melting ($E'(\tau)>0$ for $\tau\in(t_b,t_c)$) from the phases when it is growing ($E'(\tau)<0$ for $\tau\in(t_c,t_b+1)$), since $-E$ is proportional to ice thickness when $E<0$. It follows that the minimum ice thickness over the yearly cycle occurs at time $\tau=t_c$.

Letting $E_b\equiv E(t_b)$ be the initial condition and $E_c\equiv E(t_c),$ we find the period-one Poincar\'e return map, given implicitly by
\begin{equation}\label{eq:cov-poincare}
	E(t_b+1) - \frac{E(t_b+1)^2}{2\zeta} = E_b - \frac{E_c^2}{2\zeta} + \langle F_- \rangle,
\end{equation}
where
\begin{align}
	E_c &= E_b + I_-(t_b,t_c), \label{eq:Ec}\\
	\langle F_- \rangle &\equiv \int_0^1 F_-(\tau)d\tau,  \\
	I_-(t_b,t_c) &\equiv \int_{t_b}^{t_c} F_-(\tau)d\tau \label{eq:I-}.
\end{align}
The initial condition for a periodic solution, $E(t_b)=E(t_b+1)\equiv E_b^*$, is then
\begin{equation}\label{eq:Eb*}
	E_b^*=\frac{\zeta\langle F_-\rangle}{I_-(t_b,t_c)}-\frac{I_-(t_b,t_c)}{2}.
\end{equation}

We obtain the condition for a perennially ice-covered periodic state to exist by insisting that $E<0$ for all $\tau\in [t_b,t_b+1)$. Since the minimum ice thickness occurs at $\tau=t_c$, we obtain the existence condition as 
\begin{equation}\label{eq:ice-cov-exist}
	E_c^*=\langle F_-\rangle+\frac{I_-(t_b,t_c)^2}{2\zeta}<0,
\end{equation}
where $E_c^*$ is the value of $E(t_c)$ for the periodic solution beginning at $E_b^*.$ We can confirm that this solution is stable by computing the Floquet multiplier from the map (\ref{eq:cov-poincare}). We find
\begin{equation}
	\frac{dE(t_b+1)}{dE_b}\bigg\vert_{E_b=E_b^*}=\frac{\zeta-E_c^*}{\zeta-E_b^*}.
\end{equation}
Since $E_b^*<E_c^*<0$ and $\zeta$ is positive, the Floquet multiplier must be in the interval $(0,1)$, and hence the solution is stable when it exists.

Again, we have avoided the sliding intervals of the Filippov system by insisting that $E_{max}=E_c^*<0$. However, the existence boundary in parameter space corresponds to where $E(t_c)\rightarrow 0$, which is the unique entry point to $S_2$, the autumn sliding interval in Figure \ref{fig:fplusfminus}(c). Thus the bifurcation that leads to the creation (or annihilation) of the perennially ice-covered solutions also corresponds to a grazing-sliding bifurcation.

\subsection{Seasonally Ice-Free Periodic Solutions}
\label{sec:seasicefree}

For this state there are three thermodynamic phases, each governed by a different differential equation, that are fully summarized by Table~\ref{tab:seasonaldiag}. We take as a starting time for the periodic solutions $t_b$, when the ablation phase begins; our initial condition is $E_{b}\equiv E(t_b)$. We also introduce a ``melting time" $t_m$ when the ice-free phase starts, and a ``freezing time" $t_f$ when the ice-free phase ends. The time $t_b$ is defined implicitly as in (\ref{eq:tbtcdef}) (see Figure~\ref{fig:fplusfminus}(a)), while $t_m$ and $t_f$ are defined by 
\begin{align}\label{eq:tmtfdefinitions}
	\begin{aligned}
		E(t_m)&=0,\quad \frac{dE}{d\tau}(t_m)>0,\\
		E(t_f)&=0,\quad \frac{dE}{d\tau}(t_f)<0,
	\end{aligned}
\end{align}
and are not known a priori. These times are by-products of our analysis, and they determine the fraction of the year that this seasonal state is ice-free, $t_f-t_m$. Note that our requirements that $dE/d\tau>0$ at $\tau=t_m$ and $dE/d\tau<0$ at $\tau=t_f$ ensure that this solution does not enter either sliding interval at the $E=0$ boundary. This follows from the property that the sliding intervals are repelling, so they may only be entered if $dE/d\tau=0$ at their left boundaries (see Figure~\ref{fig:fplusfminus}).

\begin{table}
\centering
\setlength{\tabcolsep}{1em} 
{\renewcommand{\arraystretch}{1.5}
\begin{tabular}{|l|c|c|c|c|c|}
	\hline 
	$\vphantom{\left(\frac{E}{E-\zeta}\right)}$\textbf{Time} & $t_{b}\leq \tau<t_{m}$ &  & $t_{m}<\tau<t_{f}$ &  & $t_{f}<\tau< t_{b}+1$\tabularnewline
	\hline 
	$\vphantom{\left(\frac{\frac{1}{2}}{E-\zeta}\right)}$\textbf{Equation} & $\frac{dE}{d\tau}=F_{-}$ & \multirow{3}{*}{\begin{sideways}
  $ E\left(t_{m}\right)=0$\,\,\,\,
	\end{sideways}} & $\frac{dE}{d\tau}=F_{+}-BE$ & \multirow{3}{*}{\begin{sideways}
	$E\left(t_{f}\right)=0$\,\,\,\,
	\end{sideways}} & $\frac{dE}{d\tau}=\frac{\zeta F_-}{\zeta - E}$\tabularnewline
	\cline{1-2} \cline{4-4} \cline{6-6} 
	$\vphantom{\left(\frac{E}{E-\zeta}\right)}$\textbf{State} & $T=0,\, E<0$ &  & $E=T>0$ &  & $E,\, T<0$\tabularnewline
	\cline{1-2} \cline{4-4} \cline{6-6} 
	$\vphantom{\left(\frac{E}{E-\zeta}\right)}$\textbf{Ice Phase} & Ice melting &  & No ice &  & Ice growing\tabularnewline
	\hline
\end{tabular}
}
\caption{The three thermodynamic phases of a seasonally ice-free periodic solution, and the piecewise-defined differential equation that governs the solution.}
\label{tab:seasonaldiag}
\end{table}

Starting with $E_b\equiv E(t_b)$ and integrating the governing piecewise-defined differential equation through each of the phases indicated by Table \ref{tab:seasonaldiag}, we find that
\begin{align}\label{eq:seasonaltrajectory}
	\begin{aligned}
		E_b&=-\int_{t_b}^{t_m}F_-(\tau)d\tau \quad\,\,\left(\equiv -I_-(t_b,t_m)\right), \\
		0&=\int_{t_m}^{t_f} e^{Bt} F_+(\tau)d\tau \quad\left(\equiv I_+(t_m,t_f)\right), \\
		E(t_b+1)-\frac{E(t_b+1)^2}{2\zeta}&=\int_{t_f}^{t_b+1} F_-(\tau)d\tau \quad\,\,\left(\equiv I_-(t_f,t_b+1)\right).
	\end{aligned}
\end{align}
From these we determine that a periodic solution $E_b=E(t_b+1)\equiv E_b^*<0$ exists, provided we can find times $t_m$ and $t_f$ which solve the following implicit equations:
\begin{eqnarray}
	0&=&I_+(t_m,t_f), \label{eq:tmtfcondition1}\\
	0&=&I_-(t_b,t_m)+\frac{I_-(t_b,t_m)^2}{2\zeta}+I_-(t_f,t_b+1). \label{eq:tmtfcondition2}
\end{eqnarray}
Moreover, any solution $(t_m,t_f)$ to this pair of equations ensures existence of a periodic seasonal state only if it lies in an appropriate domain determined by the following additional constraints, related to the timeline of Table \ref{tab:seasonaldiag}:
\begin{align}\label{eq:tmtfinequalities}
	\begin{aligned}
	&t_b<t_m<t_f<t_b+1,\\
	&F_-(\tau)>0 \quad \forall\,\tau\in(t_b,t_m],\\
	&F_+(t_f)<0,
	\end{aligned}
\end{align}
where $F_-(t_b)=F_-(t_b+1)=0$. Other conditions suggested by Table~\ref{tab:seasonaldiag}, that $F_-(\tau)<0$ for $\tau\in [t_f,t_b+1)$ and $F_+(t_m)>0,$ are guaranteed by those included in~(\ref{eq:tmtfinequalities}) (see Figure \ref{fig:fplusfminus}(a)). Once the times $t_m$ and $t_f$ satisfying (\ref{eq:tmtfcondition1})--(\ref{eq:tmtfinequalities}) are determined, the initial condition $E_b^*$ for the periodic seasonally ice-free state is then given by the explicit expression in (\ref{eq:seasonaltrajectory}). 

Unlike the perennially ice-free and perennially ice-covered states, the seasonally ice-free periodic solutions are not necessarily unique~\cite{Abbot11,Eisenman09,Eisenman12}, and they need not be stable. We compute the Floquet multiplier $dE(t_b+1)/dE_b,$ evaluated at $E_b=E_b^*$, by implicitly differentiating (\ref{eq:seasonaltrajectory}) with respect to $E_b$, and find
\begin{equation}\label{eq:fmseasonal}
	\left.\frac{dE(t_b+1)}{dE_b}\right|_{E_b=E_b^*} = e^{-B(t_f-t_m)}\ \left( \frac{ \zeta}{\zeta-E_b^*}\right)\ \frac{F_-(t_f)}{ F_+(t_f)}\  \frac{ F_+(t_m) }{ F_-(t_m)}.
\end{equation}
Saddle-node bifurcations of the seasonally ice-free periodic solutions occur whenever
\begin{equation}\label{eq:sn} 
	e^{-B(t_f-t_m)}\ \left( \frac{ \zeta}{\zeta-E_b^*}\right)\ \frac{F_-(t_f)}{ F_+(t_f)}\  \frac{ F_+(t_m) }{ F_-(t_m)}=+1.
\end{equation}
Note that if $t_f\rightarrow t_a$ or $t_f \rightarrow t_d,$ or if $t_m\rightarrow t_b$ or $t_m \rightarrow t_c,$ then the Floquet multiplier would diverge since $F_+(t_a)=F_+(t_d)=0$ and $F_-(t_b)=F_-(t_c)=0$ (see Figure \ref{fig:fplusfminus}). This is relevant, as we will see further in this section and in Section \ref{sec:case-study}, near certain existence boundaries for the seasonal solutions, and causes them to come into existence as highly unstable solutions when they arise through grazing-sliding bifurcations. Additionally, note that the dimensionality of this problem precludes the possibility of a negative Floquet multiplier and period-doubling bifurcations, at least insofar as any such solutions avoid the sliding intervals, where the uniqueness of solutions breaks down. 

Finally, we may determine the seasonal state existence region in parameter space by considering all possible mechanisms for creation of seasonal ice-free solutions:
\begin{enumerate}[nolistsep]
\item{} Saddle-node bifurcations for which (\ref{eq:sn}) is met.
\item{} Solutions $(t_m,t_f)$ to (\ref{eq:tmtfcondition1}),(\ref{eq:tmtfcondition2}) entering through the boundaries of ($t_m,t_f$)-space associated with the inequalities (\ref{eq:tmtfinequalities}), which are $t_m=t_b$, $t_m=t_c$, $t_f=t_d$, and $t_f=t_a+1$. These follow from the conditions $F_-(t_m)>0$ and $F_+(t_f)<0$ in (\ref{eq:tmtfinequalities}). Below we argue that $(t_m,t_f)$ solutions to (\ref{eq:tmtfcondition1}),(\ref{eq:tmtfcondition2}) may in fact only enter or leave the ($t_m,t_f$) region of existence at the boundaries $t_m=t_c$ and $t_f=t_a+1$, which correspond to solutions that enter a sliding interval and hence are associated with grazing-sliding bifurcations.
\end{enumerate}

To eliminate the possibility of $(t_m,t_f)$ solutions to (\ref{eq:tmtfcondition1}),(\ref{eq:tmtfcondition2}) emerging from the boundaries $t_m=t_b$ and $t_f=t_d$, we consider the relationship between these times and the sliding intervals of the discontinuity boundary, as illustrated by Figure~\ref{fig:fplusfminus}. The only way a solution can have a zero at either $t_b$ or $t_d$ is if it enters the sliding interval at $t_a$ or $t_c,$ respectively. However, since the vector field at $\tau=t_a$ points downward for $E\rightarrow 0^-$(see Figure \ref{fig:fplusfminus}(b)), $t_b$ cannot correspond to a melting time. By a similar argument, $t_d$ cannot be a freezing time (see Figure \ref{fig:fplusfminus}(c)). Hence $(t_m,t_f)$ solutions do not emerge from the boundaries $t_m=t_b$ and $t_f=t_d$. The possibility of $(t_m,t_f)$ solutions emerging from the boundaries $t_m=t_c$ and $t_f=t_a+1$ is confirmed through the considerations above, namely, that $t_a$ could correspond to a freezing time and $t_c$ could correspond to a melting time. Note that these limiting values are associated with the onset of the springtime and autumn sliding intervals, respectively. Thus, $t_a$ can be considered the last possible freezing time and $t_c$ the last possible melting time.

To summarize, in addition to saddle-node bifurcations, seasonal solutions may come into existence in parameter space when we can find $(t_m,t_f)$ which satisfy
\begin{align}\label{eq:tc-tf-boundary}
	\begin{aligned}
		0&=I_+(t_c,t_f), \\
		0&=I_-(t_b,t_c)+\frac{I_-(t_b,t_c)^2}{2\zeta}+I_-(t_f,t_b+1) 
	\end{aligned}
\end{align}
at the $t_m=t_c$ boundary, or 
\begin{align}\label{eq:tm-tap1-boundary}
	\begin{aligned}
		0&=I_+(t_m,t_a+1), \\
		0&=I_-(t_b,t_m)+\frac{I_-(t_b,t_m)^2}{2\zeta}+I_-(t_a,t_b) 
	\end{aligned}
\end{align}
at the boundary where $t_f=t_a+1.$ These solutions come into existence via grazing-sliding bifurcations and are similar to the generic case $TC_1$ in \cite{Kuznetsov03}, although with repelling sliding instead of attracting sliding.

\section{Case Study Bifurcation Analysis}
\label{sec:case-study}

In this section we develop a case study comparing bifurcation results from the nondimensionalized albedo-smoothed model given by Eisenman~\cite{Eisenman12} with those associated with its Filippov limit achieved by taking the parameter $\Delta E \rightarrow 0.$ We demonstrate the similarity of analytical results in the Filippov limit to numerical results of the corresponding smoothed system from \cite{Eisenman12} through selected bifurcation diagrams.
The section concludes with an exploration of how varying the relative widths of the two sliding intervals in the Filippov system affects the bifurcation diagram.

\subsection{Characterization of Bifurcations}
\label{sec:bif-diags}

\begin{figure}[!t]
	\centering
	\includegraphics[scale=1]{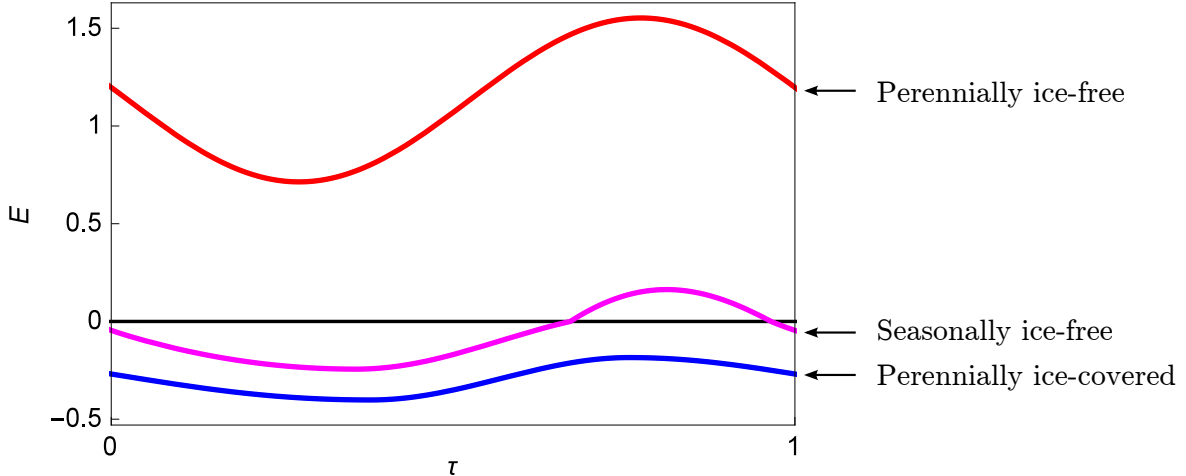}
	\caption{Examples of stable periodic solutions of (\ref{eq:system-summary}), where $F_\pm(\tau)$ are given by (\ref{eq:PWSFE12}), for the default parameters $S_a=1.5,$ $L_a=0.73,$ $\phi=0.15$, $B=0.45,$ $\zeta=0.12,$ and $\Delta\alpha=0.43,$ as in Table 1 of \cite{Eisenman12}. We set $L_m=0.92$ for the perennially and seasonally ice-free trajectories and $L_m=1.1$ for the perennially ice-covered trajectory.}
	\label{fig:solution-trajectories}
\end{figure}

The simplified Eisenman model is given by (\ref{eq:dEdt}), (\ref{eq:TE12}), with sinusoidal approximations to the periodic forcing functions~\cite{Eisenman12}:
\begin{equation}\label{eq:FE12}
	F(\tau,E) = (1 + \Delta\alpha\tanh(E/\Delta E))(1-S_a \cos(2\pi\tau)) - (L_m + L_a \cos(2\pi(\tau-\phi))).
\end{equation}
Here, $\Delta\alpha>0$ is proportional to the jump in albedo between sea ice ($E<0$) and open ocean ($E \ge 0$) as defined in~(\ref{eq:Da}), $\Delta E$ is the energy scale for smoothing, $S_a>0$ captures the seasonal variation of incoming shortwave radiation, $L_a>0$ measures the seasonal variation of longwave radiation, $\phi\in [-\frac{1}{2},\frac{1}{2})$ controls the phase shift between the peaks in incoming shortwave and outgoing longwave radiation, and $L_m>0$ is the mean of the outgoing longwave radiation when $T=0$. Note that in \cite{Eisenman12}, $\phi\in[0,1);$ we have chosen to center the interval at $\phi=0.$ In the Filippov limit $F(\tau,E)$ takes the form
\begin{equation}\label{eq:PWSFE12}
	F(\tau,E) = 
	\begin{cases}
		F_+(\tau) \equiv (1 + \Delta\alpha)(1-S_a \cos(2 \pi \tau)) - (L_m + L_a \cos(2 \pi (\tau-\phi))), & E>0, \\
		F_-(\tau) \equiv (1 - \Delta\alpha)(1-S_a \cos(2 \pi \tau)) - (L_m + L_a \cos(2 \pi (\tau-\phi))), & E<0.
	\end{cases}
\end{equation}
Figure \ref{fig:solution-trajectories} shows examples of the three types of periodic solutions constructed from the analytical results of Section \ref{sec:solutions}, using the case study model (\ref{eq:system-summary}),(\ref{eq:PWSFE12}) for parameter sets where the solutions exist stably.

In our bifurcation analysis, as in previous studies \cite{Abbot11,Eisenman12,Moon11,Moon12} of the model introduced by Eisenman and Wettlaufer \cite{Eisenman09}, we consider a decrease in the mean value of the outgoing longwave radiation (for $T=0$) to serve as a suitable proxy for increases in greenhouse gas emissions. For example, the mean value of $F_0(t)$ in (\ref{eq:EW09}), which is denoted $L_m$ after nondimensionalization in (\ref{eq:PWSFE12}), is expected to decrease approximately logarithmically with increases in atmospheric CO$_{\text{2}}$ concentration \cite{Pierrehumbert11}. Thus, increases in greenhouse gases have the effect of decreasing the widths of the sliding intervals in Figure \ref{fig:fplusfminus}. The default present-day value of $L_m$ is given as 1.25 in \cite{Eisenman12} and is estimated based on observations \cite{Eisenman09,Eisenman12}.

The Eisenman model (\ref{eq:dEdt}),(\ref{eq:TE12}),(\ref{eq:FE12}) captures the same bifurcation behavior as the model in EW09, as well as the results of several other single-column energy balance and global climate models \cite{Eisenman12}. Through extensive numerical simulations for various parameter sets, Eisenman identified four bifurcation scenarios, representing the range of results from these previous models, which also occur in his model. For the purposes of our analysis, this model has the benefit that in the Filippov limit the integrals associated with solutions, derived in Section \ref{sec:solutions}, can be evaluated in closed form. On the other hand, for some parameter sets the model violates the assumption that $F_+(\tau)\geq F_-(\tau),$ for all $\tau,$ made throughout previous sections. In this analysis we will not consider parameter sets for which this violation leads to attracting sliding intervals.

\begin{figure}[!t]
	\centering
	\begin{tabular}{ccc}
	\multicolumn{2}{c}{\subfloat[]{\includegraphics[scale=1]{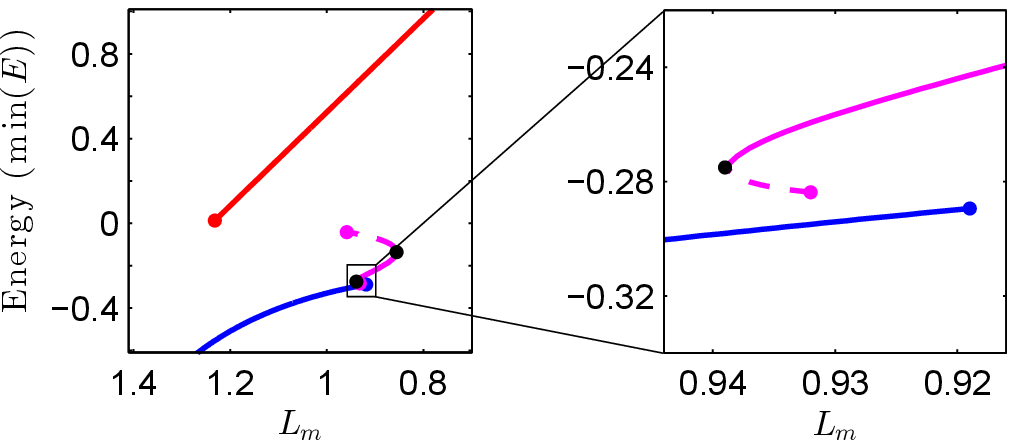}}} & 
	\subfloat{\includegraphics[scale=1, trim=0 -60 0 0, clip=true]{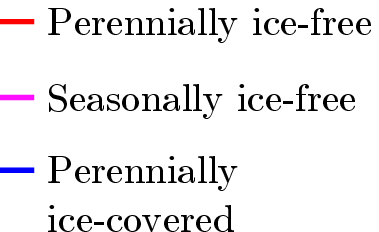}} \\
	\renewcommand{\thesubfigure}{b}
	\subfloat[]{\includegraphics[scale=1]{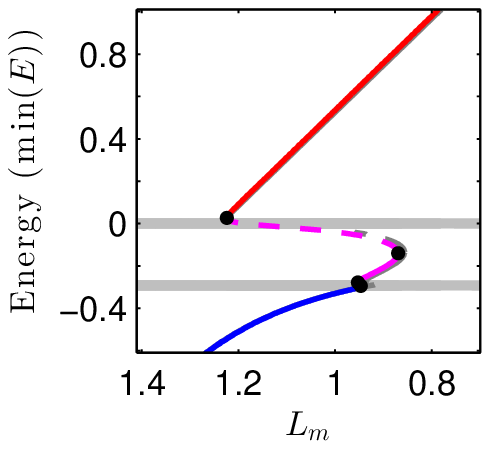}} &
	\renewcommand{\thesubfigure}{c}
	\subfloat[]{\includegraphics[scale=1]{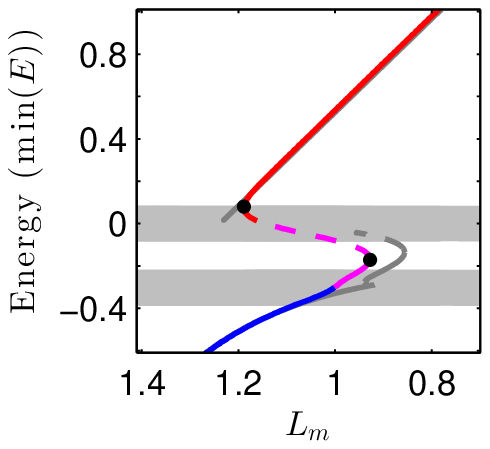}} &
	 \\
	\end{tabular}
\caption{Bifurcation diagrams for (\ref{eq:dEdt}),(\ref{eq:TE12}), where $F(\tau,E)$ is given by (\ref{eq:FE12}), with (a)~$\Delta E = 0$, (b)~$\Delta E=0.02$ (colored curves) and $\Delta E = 0$ (gray curves), and (c) $\Delta E=0.08$ (colored curves) and $\Delta E = 0$ (gray curves), for the default parameters, as in Figure \ref{fig:solution-trajectories}. In (a) the diagram on the right is a blow-up of the boxed region in the left diagram. Black (colored) points denote saddle-node (grazing-sliding) bifurcations. Solid (dashed) curves represent stable (unstable) solutions. $L_m$ decreases with increasing greenhouse gases. Horizontal shading around each transition in (b)--(c) indicates the amount of the smoothing (shading width is $2\Delta E$).}
\label{fig:bif-diag-comparison}
\end{figure}

An example bifurcation diagram of the Filippov system is shown in Figure \ref{fig:bif-diag-comparison}(a). Two bifurcation diagrams of the albedo-smoothed system are given for comparison, in (b) with a small amount of smoothing ($\Delta E=0.02$), and in (c) with the ``default'' amount of smoothing used in the analysis of \cite{Eisenman12} ($\Delta E=0.08$). We make the following remarks about Figure~\ref{fig:bif-diag-comparison}:
\begin{enumerate}[nolistsep]
\item Since $L_m$ decreases with increased greenhouse gases, we have chosen to show $L_m$ decreasing from left to right in the bifurcation diagram (so that greenhouse gases increase from left to right). Solutions progress from the current state of a perennially ice-covered Arctic, through a seasonally ice-covered state, to a state in which the Arctic is perennially ice-free at high greenhouse gas levels. The transition from seasonal ice to a perennially ice-free state results from a saddle-node bifurcation on the seasonal branch.
\item Notice that we have chosen to represent the periodic solution branches in the bifurcation diagrams using $\min(E)$ instead of values chosen from a single Poincar\'e map.
\item We calculate the bifurcation diagram numerically for the albedo-smoothed system in (b) and (c) by integrating a range of initial conditions for one year to generate a Poincar\'e map $E_{n+1} = F(E_n)$, as in \cite{Eisenman09,Eisenman12}; each zero of the function $G(E)=E-F(E)$ gives an initial value for a periodic state, and these initial conditions are used to compute $\min(E)$. For the bifurcation diagram for the Filippov system, we use the analytical solutions calculated in Section \ref{sec:solutions}. For the seasonal state, the melting and freezing times $t_m$ and $t_f$, defined implicitly by (\ref{eq:tmtfcondition1}),(\ref{eq:tmtfcondition2}), are found numerically.
\item There are two gaps in the Filippov bifurcation diagram in (a); these gaps originate when solutions limit to the sliding intervals, as described in Section \ref{sec:solutions}. (Recall that solutions that involve sliding are nonunique; we do not depict their corresponding solution branches.) The termination points of the solution branches in (a) occur where solutions limit to the sliding intervals; they correspond to grazing-sliding bifurcations.
\item The seasonal state of the Filippov system always comes into existence through grazing-sliding as an unstable branch. This follows from the Floquet multiplier (\ref{eq:fmseasonal}), which diverges at the existence boundaries associated with solutions limiting to the sliding intervals.
\item In (a) and (b), notice that there are two hysteresis loops: a tiny one that involves the transition between the perennially ice-covered solution branch and the seasonally ice-free branch, and a large one that involves the perennially ice-free branch. The smaller hysteresis loop is shown in detail in the blown-up region in (a). This hysteresis loop is easily smoothed out by increasing the amount of smoothing in the model and is not present in (c). The absence of this bifurcation associated with the loss of summer sea ice was a major result of EW09.
\item Most qualitative aspects of the bifurcation diagram of the albedo-smoothed system are preserved in the Filippov limit; notice that this resemblance increases as $\Delta E$ decreases, which is illustrated in the difference between (b) and (c). (See Appendix \ref{app:comparison} for a more thorough comparison of bifurcations of the Filippov system and the albedo-smoothed system.) 
\end{enumerate}

\subsection{Impact of sliding intervals}
\label{sec:simplified-model}

\begin{figure}[!t]
	\centering
	\begin{tabular}{cc}
	\subfloat[]{\includegraphics[scale=1]{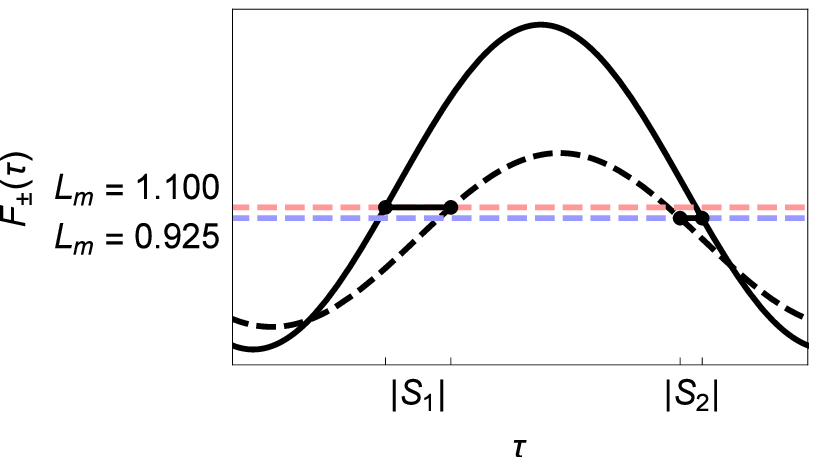}} &
	\subfloat[]{\includegraphics[scale=1]{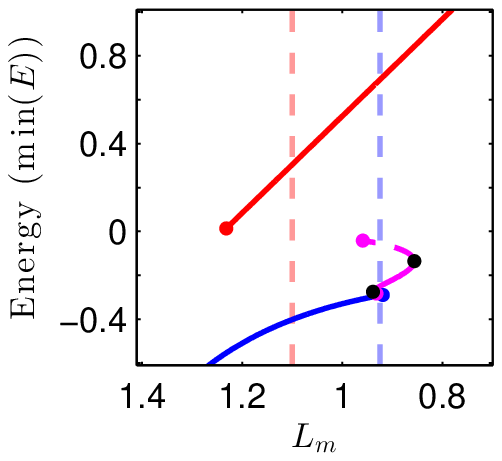}} \\
	\subfloat[]{\includegraphics[scale=1]{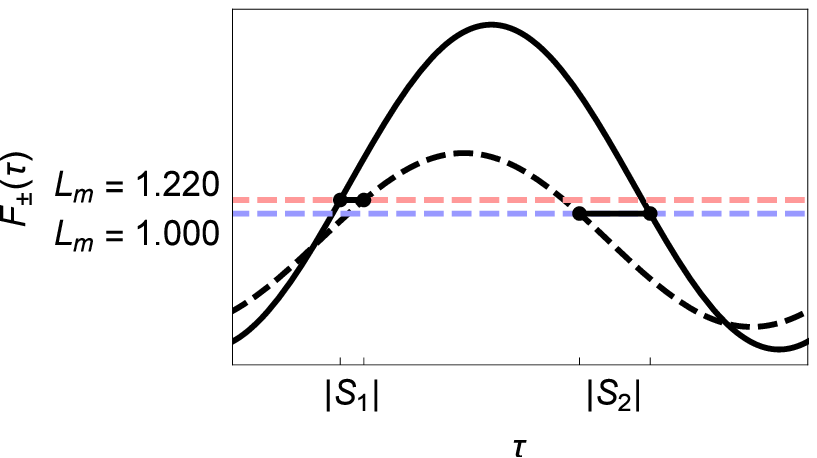}} &
	\subfloat[]{\includegraphics[scale=1]{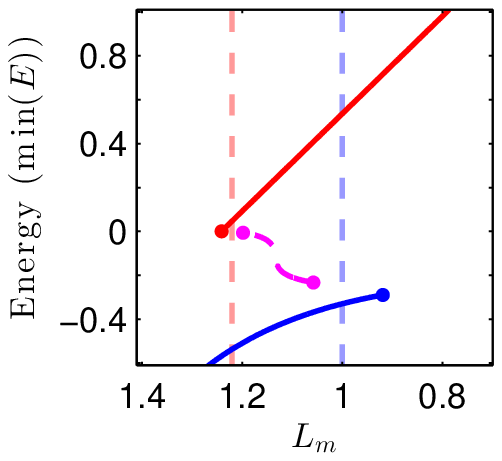}} \\
	\end{tabular}
	\caption{Illustration of the relationship between the sliding interval widths and the bifurcation diagram gaps. (a),(c) show $F_+(\tau)$ (solid curve) and $F_-(\tau)$ (dashed curve) for the Filippov system (\ref{eq:system-summary}),(\ref{eq:PWSFE12}); $F_\pm=0$ lines for different $L_m$ values are indicated by horizontal dashed lines. (b),(d) show corresponding bifurcation diagrams. Parameters of (\ref{eq:PWSFE12}) are as in Figure \ref{fig:solution-trajectories}, except in (c),(d) we set $S_a=1.58,$ $L_a=0.85,$ and $\phi=-0.20$. The sliding interval corresponding to each bifurcation diagram gap is indicated by a solid line in (a),(c) for the $L_m$ values indicated by vertical lines in (b),(d).}
	\label{fig:sliding-region-intuition}
\end{figure}

In this section we explore how the widths of the sliding intervals impact the size of the hysteresis loops in the corresponding albedo-smoothed system. We start with the observation that, roughly speaking, a wide sliding interval leads to a wide gap in the bifurcation diagram, whereas a narrow sliding interval creates a narrow gap. We find that for some parameter sets, narrow gaps may be associated with small hysteresis loops that are easily removed when smoothing is introduced. Figure \ref{fig:sliding-region-intuition} provides two side-by-side examples demonstrating how the widths of the sliding intervals impact the sizes of the bifurcation diagram gaps. We use these examples as motivation for considering how the relationship between $F_+$ and $F_-$ impacts the bifurcation structure.

We wish to explore how changing the relative widths of the spring and autumn sliding intervals affects the bifurcation diagram. For this purpose we find it more convenient to write $F_\pm$ given by (\ref{eq:PWSFE12}) in the following standard form:
\begin{equation}\label{eq:simple-FE12}
	F_\pm(\tau) = \overline{F_{\pm}} + \widetilde{F_\pm} \cos(2\pi \tau - \psi_{\pm}),
\end{equation}
where
\begin{eqnarray}
	\overline{F_\pm} &=& 1-L_m\pm\Delta\alpha, \nonumber\\
	\widetilde{F_\pm} \cos(\psi_\pm) &=& -(1\pm\Delta\alpha)S_a-L_a\cos(2\pi\phi), \\
	\widetilde{F_+}\sin(\psi_+) &=& \widetilde{F_-}\sin(\psi_-) = -L_a\sin(2\pi\phi), \nonumber
\end{eqnarray}
from which it follows that
\begin{align}
	\begin{aligned}
		\widetilde{F_\pm} &= \sqrt{(1\pm\Delta\alpha)^2 S_a^2 + 2(1\pm\Delta\alpha)S_a L_a \cos(2\pi\phi) + L_a^2}  \\
		\Delta\psi &\equiv \psi_+ - \psi_- \\
			&= -\mbox{sgn}(\phi)\arccos\left( \frac{ \left( 1-\Delta\alpha^2 \right)S_a^2 + 2S_a L_a \cos(2\pi\phi)+L_a^2 }{ \widetilde{F_+}\widetilde{F_-} } \right),
	\end{aligned}
\end{align}
where the factor $-\mbox{sgn}(\phi)$ ensures that we have inverted the cosine function correctly; note that $\Delta\psi=0$ if $\phi=0$.

The form of $F_\pm$ in (\ref{eq:simple-FE12}) has the property that each parameter change causes a transparent transformation of the periodic function. It suggests a way to sweep through the parameter space of the model that is different from the parameter studies performed by Eisenman in \cite{Eisenman12} and is more systematic for our investigation of the impact of the sliding intervals. Note that the inverse of this mapping is not always unique: each parameter set given by (\ref{eq:PWSFE12}) has a corresponding set in (\ref{eq:simple-FE12}), but parameter sets in (\ref{eq:simple-FE12}) may map to parameters in (\ref{eq:PWSFE12}) which are not permitted due to physical restrictions, or they may map to two sets of parameters. See Appendix \ref{app:inverse-mapping} for a description of the inverse mapping we use in our investigation.

\begin{figure}[!t]
	\centering
	\begin{tabular}{ccc}
	\subfloat[]{\includegraphics[scale=1]{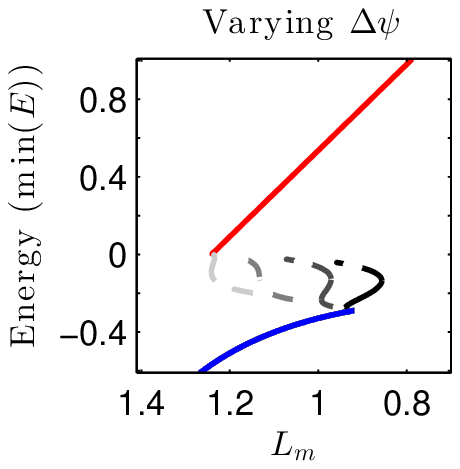}} &
	\subfloat[]{\includegraphics[scale=1]{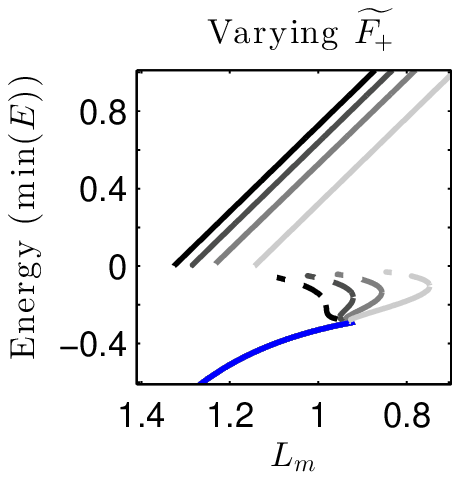}} &
  \subfloat[]{\includegraphics[scale=1]{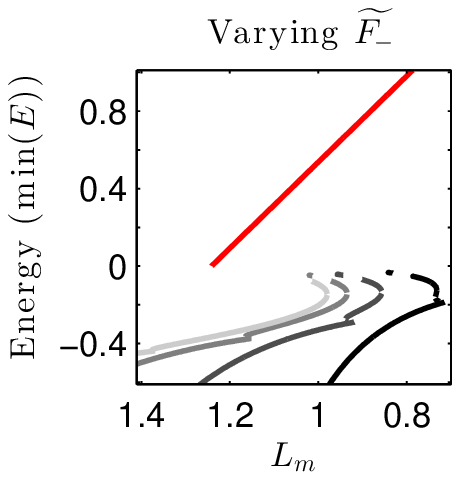}} \\
	\subfloat[]{\includegraphics[scale=1]{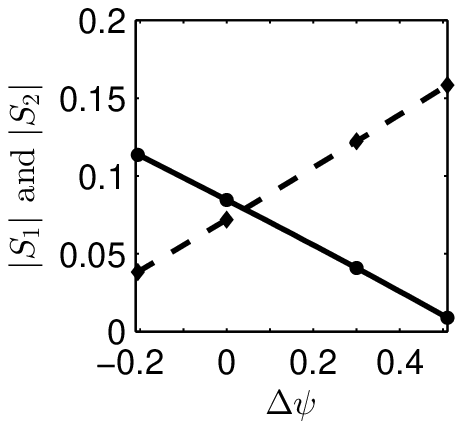}} &
	\subfloat[]{\includegraphics[scale=1]{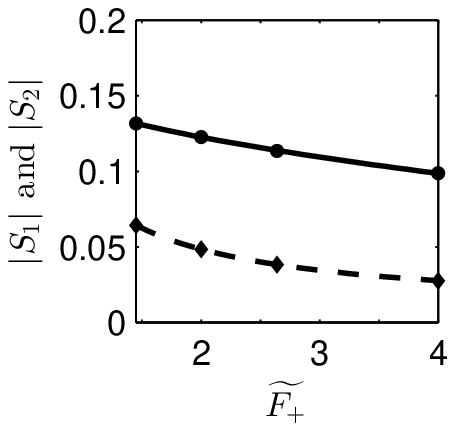}} &
  \subfloat[]{\includegraphics[scale=1]{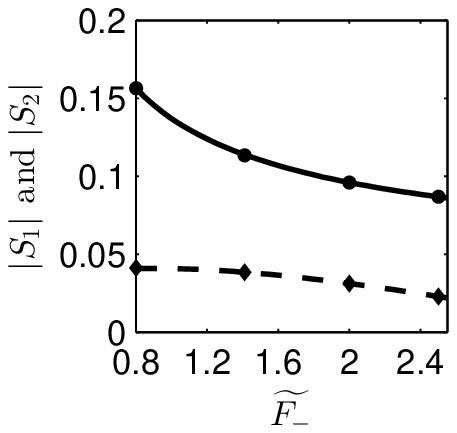}} \\
	\end{tabular}
	\caption{(a)--(c) Bifurcation diagrams of the Filippov system (\ref{eq:system-summary}),(\ref{eq:PWSFE12}), with branches which change as an additional parameter of (\ref{eq:simple-FE12}) is varied overlaid in shades of gray, corresponding to black, dark gray, medium gray, and light gray branches as each parameter is increased (unchanged branches retain their original color). Diagrams correspond to (a) $\Delta\psi = -0.21$ (default parameters), $0.00,$ $0.30,$ and $0.51$; (b) $\widetilde{F_+}=1.45,$ $2.00,$ $2.64$ (default), and $4.00$; and (c) $\widetilde{F_-}=0.80,$ $1.41$ (default), $2.00,$ and $2.50.$ Parameters used to generate each diagram for the value of $\Delta\psi$ or $\widetilde{F_\pm}$ shown are given in Table~\ref{tab:bif-pars} in Appendix~\ref{app:inverse-mapping}. Note: the black and medium gray seasonally ice-free branches in (a) are identical to the seasonal branches in Figure \ref{fig:sliding-region-intuition}(b),(d), respectively. (d)--(f) Curves showing the widths of the $S_1$ and $S_2$ sliding intervals (solid and dashed curves, respectively) as each parameter is varied. Markers highlight the sliding interval widths for the parameter values used in the corresponding figure above.}
	\label{fig:vary-DPsi-and-Fpmtilde}
\end{figure}

Now we can explore how parameters in $F_\pm$ change the bifurcation diagram. See Figure~\ref{fig:vary-DPsi-and-Fpmtilde} for several examples of how varying parameters of (\ref{eq:simple-FE12}) in turn, (a) $\Delta\psi,$ (b) $\widetilde{F_+}$, and (c)~$\widetilde{F_-},$ affect the bifurcation diagram. We make the following observations about Figure~\ref{fig:vary-DPsi-and-Fpmtilde}:
\begin{enumerate}[nolistsep]
\item The parameter $\widetilde{F_+}$ only affects the seasonal and perennially ice-free solution branches, while $\widetilde{F_-}$ affects only the seasonal and perennially ice-covered states.

\item The phase shift $\Delta\psi\equiv \psi_+ - \psi_-,$ in contrast to $\widetilde{F_\pm},$ only impacts the seasonal solution branches of the bifurcation diagram and is perhaps the most interesting to vary. For example, Figure \ref{fig:vary-DPsi-and-Fpmtilde}(a) shows that increasing $\Delta \psi$, which controls the temporal offset between the peaks of $F_+$ and $F_-,$ causes the seasonally ice-free branch to shift left to higher $L_m$ values. This causes the higher-energy gap in the bifurcation diagram to become smaller, while the lower-energy gap becomes larger.

\item The widths of the sliding intervals, as each parameter ($\Delta\psi,\widetilde{F_{+}},\widetilde{F_-}$) is varied, are shown in Figure \ref{fig:vary-DPsi-and-Fpmtilde}(d)--(f). Here, for each parameter set we use the median value of $L_m$ in each bifurcation gap to determine the width of the corresponding sliding interval. Notice that the two sliding intervals switch their relative widths for larger values of $\Delta\psi,$ so that $|S_2|>|S_1|,$ while $|S_2|<|S_1|$ for all values of $\widetilde{F_{\pm}}$ considered. As illustrated in Figure~\ref{fig:bif-diag-comparison}, this behavior is related to the bifurcation diagram gaps changing their relative widths as $\Delta\psi$ increases.
\end{enumerate}

\begin{figure}[!t]
	\centering
	\begin{tabular}{ccc}
	\subfloat[]{\includegraphics[scale=1]{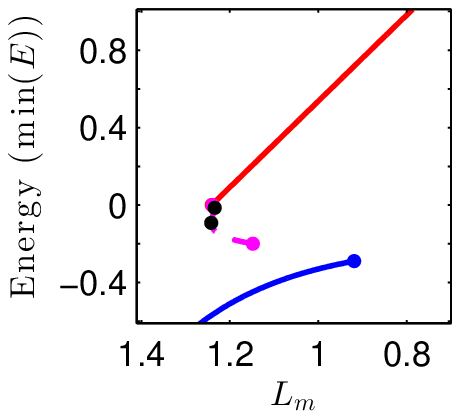}} &
	\subfloat[]{\includegraphics[scale=1]{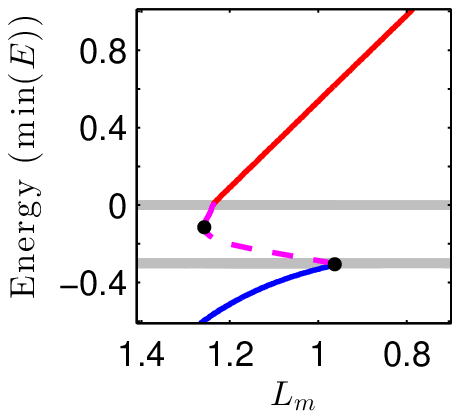}} &
  \subfloat[]{\includegraphics[scale=1]{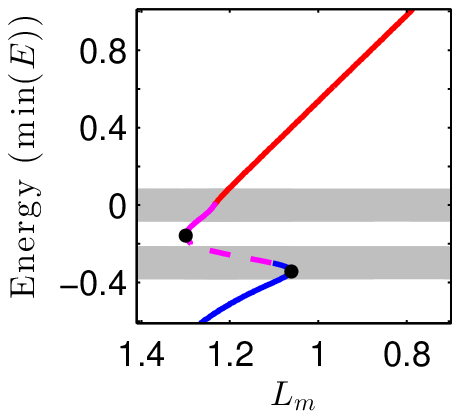}} \\
	\end{tabular}
	\caption{Bifurcation diagrams for (\ref{eq:dEdt})-(\ref{eq:TE12}), where $F(\tau,E)$ is given by (\ref{eq:FE12}), with (a)~$\Delta E = 0$, (b)~$\Delta E=0.02$, and (c) $\Delta E=0.08$, for $\Delta\psi=0.51,$ obtained for $S_a=1.82,$ $L_a=1.19,$ $\phi=-0.28$, and other default parameters as in Figure \ref{fig:solution-trajectories}.}
	\label{fig:smoothed-DPsi}
\end{figure}

As $\Delta\psi$ increases in Figure \ref{fig:vary-DPsi-and-Fpmtilde}(a), the transition involving the perennially ice-covered solution branch, with increasing greenhouse gases, becomes large. A tiny hysteresis loop may form instead, which involves the transition between the perennially ice-free branch and the seasonally ice-free branch. Figure \ref{fig:smoothed-DPsi} shows the effect of reintroducing smoothing into the albedo function on the behavior of the bifurcation diagram in this extreme case, where $\Delta\psi = 0.51$. Here, the smaller hysteresis loop that gets smoothed out is the transition from the perennially ice-free state to the seasonally ice-free state. In effect, increasing $\Delta\psi$ has moved the small hysteresis loop from the low energy, thicker maximum ice-thickness end of the seasonal branch to the high energy, thinner maximum ice-thickness end. There is now a sudden transition from the perennially ice-covered branch to the perennially ice-free branch. The smooth transition to the seasonal state with increasing greenhouse gases is completely lost, although there is now a smooth transition to a seasonal state from the ice-free state as greenhouse gases are decreased.

The example presented in Figure~\ref{fig:smoothed-DPsi} illustrates how the relative widths of the sliding intervals might impact the bifurcation diagram. It falls under Scenario 3 of Eisenman's classification, which is concerned with possible abrupt transitions under increases in greenhouse gases, starting with the present-day perennial ice state \cite{Eisenman12}. However, varying a combination of the parameters used above may lead to more variation in possible bifurcation scenarios. To this end, we consider how the size of the jump from the end of the ice-covered branch to the seasonally ice-free branch, $\Delta (\min(E)),$ varies while changing a pair of parameters. This proxy for the significance of an abrupt transition is small for the default parameters ($\Delta (\min(E))\approx 0.04$; see Figure \ref{fig:bif-diag-comparison}(a), zoomed box), but $\Delta (\min(E))$ may be increased by both decreasing $\widetilde{F_+}$ and increasing $\widetilde{F_-}$ from their default values, the combination of which is similar to a decrease in $\Delta\alpha,$ while also increasing $\Delta\psi.$

\begin{figure}[!t]
	\centering
	\begin{tabular}{ll}
		 \multirow{2}*[0.18in]{\subfloat[]{\includegraphics[scale=1]{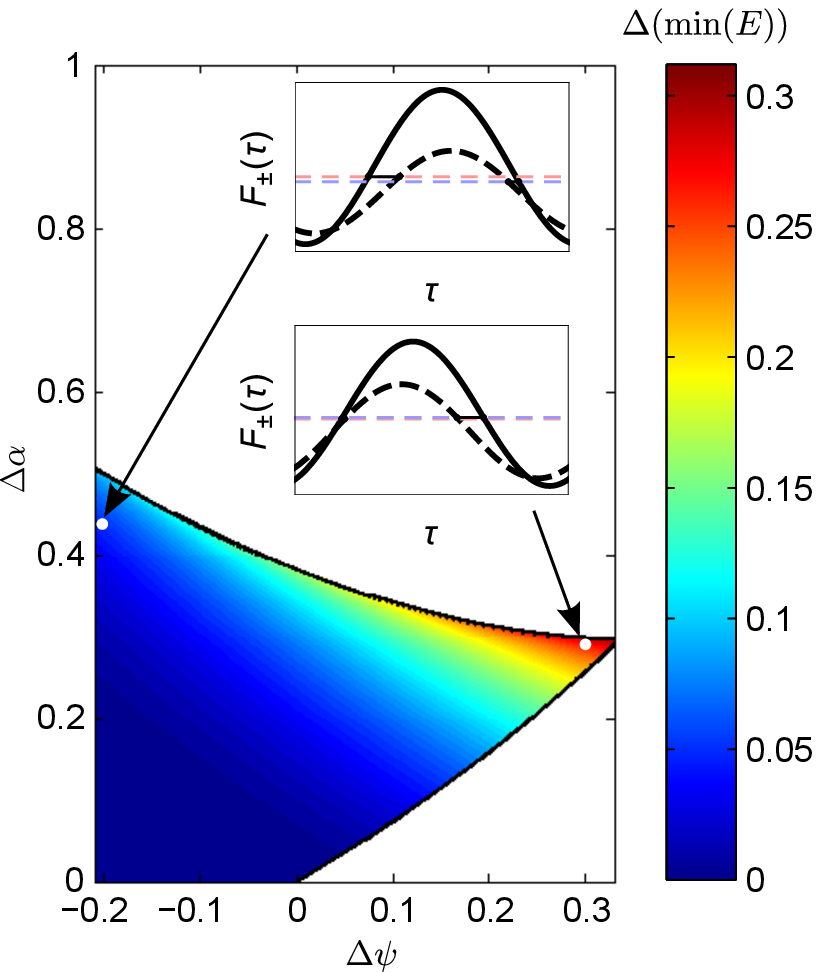}}} \\
			& \subfloat[\textcolor{white}{\,space plez}]{\includegraphics[scale=1]{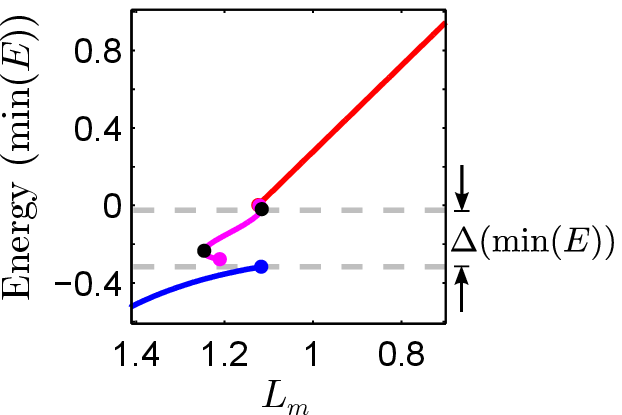}} \\
			& \vspace{0.01in} \subfloat[]{\includegraphics[scale=1]{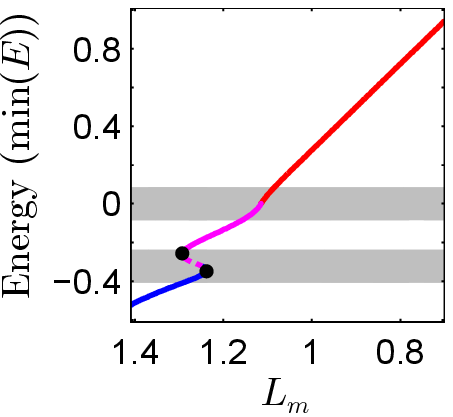}} \\
  \end{tabular}
	\caption{(a) $\Delta (\min(E))$ as a function of $\Delta\psi$ and $\Delta\alpha$; inset: $F_\pm(\tau)$ for the default parameters (top) and for the parameters used in (b),(c) (bottom). (b) and (c) show bifurcation diagrams for $(\Delta\psi,\Delta\alpha)=(0.27,0.3),$ obtained for $S_a=1.68,$ $L_a=1.03,$ and $\phi=-0.24,$ in (a) the Filippov model and (b) in the albedo-smoothed model with $\Delta E=0.08$ (cf.~the default parameters: $S_a=1.5,$ $L_a=0.73,$ $\phi=0.15$, and $\Delta\alpha=0.43$). The distance between horizontal lines in (b) gives $\Delta(\min(E))$ for $(\Delta\psi,\Delta\alpha)=(0.27,0.3).$}
	\label{fig:DE}
\end{figure}

Figure \ref{fig:DE}(a) shows a plot of $\Delta (\min(E))$ as a function of $\Delta\psi$ and $\Delta\alpha.$ A relatively large value of $\Delta (\min(E)),$ approximately 0.3, occurs at $(\Delta\psi,\Delta\alpha)\approx(0.27,0.3).$ The bifurcation diagram for these parameters is given in (b) for the Filippov model and (c) for the albedo-smoothed model. We make the following remarks about Figure \ref{fig:DE}:
\begin{enumerate}[nolistsep]
\item[1.] The horizontal dashed lines in (b) indicate where solutions jump off the ice-covered branch and hit the seasonally ice-free branch. The distance between the horizontal lines is $\Delta (\min(E)).$
\item[2.] For a large region in (a), $\Delta (\min(E))$ is similar to its value when calculated using the default parameters, where $\Delta (\min(E))\approx 0.04$. However, for high values of $\Delta\psi,$ $\Delta (\min(E))$ can become much larger, with the largest values occurring close to the upper black boundary.
\item[3.] In the white region for $\Delta\alpha>0.3$ in (a), ice-covered solutions jump directly to the ice-free branch instead of the seasonal branch, as in Figure \ref{fig:smoothed-DPsi}(a). The black boundary along the edge of this region is where the higher-energy saddle-node bifurcation on the seasonal branch coincides with the $L_m$ value of the ice-covered grazing-sliding bifurcation. The white region for lower $\Delta\alpha$ values ($\Delta\alpha<0.3$) represents where there are attracting sliding intervals; the black boundary here is where $t_c=t_d$ (see Figure~\ref{fig:fplusfminus}).
\item[4.] The default value of $\Delta\alpha=0.43$ is determined in \cite{Eisenman12} using $\alpha_i=0.68$ and $\alpha_{ml}=0.2.$ A decrease in $\Delta\alpha,$ given by (\ref{eq:Da}), to $0.3$ may be achieved, for example, using $\alpha_i=0.57$ and keeping $\alpha_{ml}=0.2$.
\item[5.] Notice that the jump to the seasonal branch in (b), $\Delta (\min(E))\approx 0.3,$ is smaller than the jump from the seasonal branch to the ice-free branch in the bifurcation diagram using the default parameter set, which is approximately $1.0$ (see Figure \ref{fig:bif-diag-comparison}(a)).
\item[6.] In (b), the jump to the seasonal state is much larger than the jump to the ice-free state. Note that there is no jump from the seasonal to the ice-free branch when the smoothing is introduced in (c).
\end{enumerate}

\section{Discussion}
\label{sec:discussion}

In our analysis we consider an energy balance model of the type introduced in EW09. The model in EW09 shows a smooth transition from a perennially ice-covered Arctic state to a seasonally ice-free state as greenhouse gases increase (see Figure \ref{fig:bif-diag-comparison}(c)). Through their analysis Eisenman and Wettlaufer showed that sea ice thermodynamics may play a key role in this smooth transition by counteracting the effects of positive ice-albedo feedback. In~\cite{Eisenman12}, Eisenman presented results of a thorough numerical exploration of bifurcation behavior in parameter space for a similar model, determining the $L_m$ values where bifurcations and transitions between solution states occur as each parameter in~(\ref{eq:FE12}) is varied in turn, while fixing the others to their default values. He identified four possible bifurcation scenarios that the model produces, and discussed their physical interpretations and connections to results of other previously studied models.

Our findings provide a case study that demonstrates how changing the relative widths of the sliding interval in a Filippov limit can affect hysteresis loops in bifurcation diagrams associated with the albedo-smoothed problem. Because we find that certain qualitative aspects of the bifurcations in the Filippov system are similar to their counterparts in the albedo-smoothed system, changes to the bifurcation diagram caused by varying the sliding interval widths are often reflected in the albedo-smoothed system as well.

One uncertainty in these models, and in similar single-column models of the Arctic, lies in the representation of the albedo transition from that of ice to open ocean. This transition is often modeled as a sigmoid (see (\ref{eq:AEW09})), where the smooth transition parameterizes the spatial extent of the patchy region of ice between the ice-covered pole and open ocean at lower latitudes, as well as the dependence of albedo on bare ice thickness. However, model results are sensitive to the amount of smoothing. Additionally, a recent study by Wagner and Eisenman~\cite{Wagner15} investigated the limitations of parameterizations like the sigmoidal albedo transition in single-column models, by combining the EW09 single-column model with a spatially varying model.

Our case study has a similar aim as previous studies of models of the type introduced in EW09: to understand how bifurcation behavior changes through parameter space. It is motivated by the relationship between the sliding interval widths and bifurcation diagram gaps, which are tied to the extent of hysteresis. We explore bifurcation behavior of the Filippov system by taking various paths through Eisenman's parameter space which were suggested by the form of the periodic forcings and how they affected the widths of the sliding intervals. This provides an alternative perspective, where the focus is on the relative phase of the seasonal forcing (for an ice-free versus an ice-covered Arctic) in the presence of sea ice self-insulation feedback, rather than just on the role of the self-insulation feedback. From this perspective, we have identified a scenario which may not be reflected in the scenarios considered in \cite{Eisenman12}, e.g.,~the example in Figure \ref{fig:DE}(c).

\begin{figure}[!t]
	\centering
	\begin{tabular}{ccc}
		\subfloat[]{\includegraphics[scale=1]{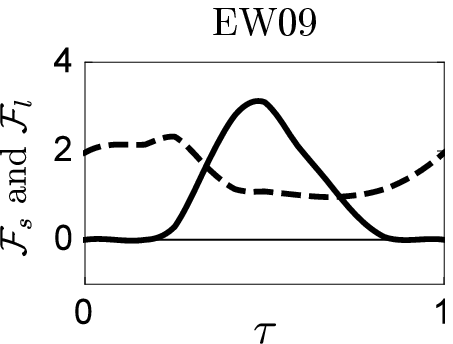}} &
		\subfloat[]{\includegraphics[scale=1]{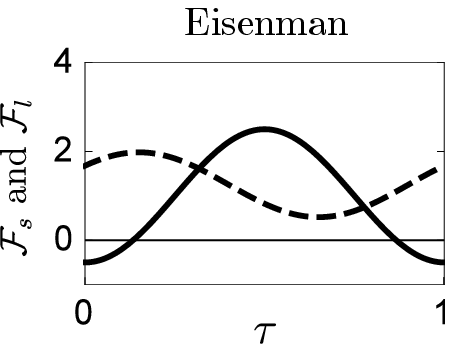}} &
		\subfloat[]{\includegraphics[scale=1]{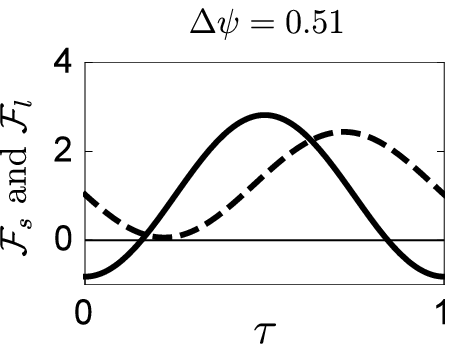}} \\
		\end{tabular}
	\caption{Incoming shortwave radiation ${\cal F}_s(\tau)$ (solid curve) and outgoing longwave radiation ${\cal F}_l(\tau)$ (dashed curve), as defined by (\ref{eq:F}) with ${\cal F}_{l,+}={\cal F}_{l,-}$. (a) From EW09 \cite{Eisenman09}, where ${\cal F}_s$ and ${\cal F}_l$ have been nondimensionalized by the annual mean incoming shortwave radiation and the mean absorbed shortwave radiation, respectively \cite{Eisenman12}. (b)~From the sinusoidal approximation of EW09 used in \cite{Eisenman12}, where ${\cal F}_s(\tau) \equiv 1-S_a \cos\left( 2\pi\tau\right)$ and ${\cal F}_{l}(\tau)\equiv L_m+L_a\cos \left(2\pi\left(\tau -\phi\right)\right)$ from (\ref{eq:PWSFE12}), using the default parameters, as in Figure \ref{fig:solution-trajectories}. (c) From (\ref{eq:PWSFE12}), using parameters, as in Figure~\ref{fig:smoothed-DPsi}, that lead to $\Delta\psi=0.51.$}
	\label{fig:Fs-and-Fl}
\end{figure}

However, the relevancy of the example bifurcation diagrams given in Section \ref{sec:simplified-model} must be considered. The parameters used in Figures \ref{fig:smoothed-DPsi} and \ref{fig:DE} are within the range considered in the parameter studies performed in \cite{Eisenman12}. The largest difference in parameter values between the default set and those used in Figures \ref{fig:smoothed-DPsi} and \ref{fig:DE} is the phase shift, $\phi$. The phase shift is the difference between when the maximum incoming shortwave radiation (${\cal F}_s$ in (2.6)) and minimum outgoing longwave radiation occur (${\cal F}_l$ in (2.6)). (Note that the \textit{minimum} of ${\cal F}_l,$ a loss term, corresponds to the \textit{maximum} of the longwave forcing.) See Figure \ref{fig:Fs-and-Fl}(a) for an example of ${\cal F}_{s,l}$ using the seasonal forcing from the model in EW09 \cite{Eisenman09}, nondimensionalized as in Eisenman's analysis \cite{Eisenman12}. Figure \ref{fig:Fs-and-Fl}(b) shows ${\cal F}_{s,l}$ from (\ref{eq:PWSFE12}) using the default parameters, as in Figure \ref{fig:solution-trajectories} ($\phi=0.15$). Note the close correspondence between the forms of ${\cal F}_{s,l}$ in (a) and (b); this is a result of using a sinusoidal approximation of the observation-derived forcing ${\cal F}_{s,l}$ in (a) to determine the form and default parameters of ${\cal F}_{s,l}$ in (b) \cite{Eisenman12}. Figure \ref{fig:Fs-and-Fl}(c) also shows ${\cal F}_{s,l}$ from (\ref{eq:PWSFE12}), but using the parameter set corresponding to $\Delta\psi=0.51$ instead, as in Figure~\ref{fig:smoothed-DPsi}. For $\phi=-0.28,$ the phase shift in (c) (as with $\phi=-0.24,$ the phase shift in Figure~\ref{fig:DE}), ${\cal F}_l$ reaches a minimum before the summer solstice, or alternatively, it is much delayed from the summer (due to the periodicity of the forcing). 

Physically, the seasonally varying component of the longwave forcing in Figure \ref{fig:Fs-and-Fl}(c) might occur in a different climate than the present-day (perennially ice-covered) state represented by the default parameter sets in both the model in EW09 \cite{Eisenman09} and the subsequent model by Eisenman \cite{Eisenman12}. There are two main contributions to the seasonal variation in the outgoing longwave radiation modeled in EW09: one contribution captures any seasonal variation in Arctic cloudiness, and the other contribution arises from the seasonal variation of atmospheric heat transport into the Arctic from lower latitudes \cite{Eisenman09}. In an alternative scenario, with perennially ice-free conditions, it has been suggested that a heavy cloud cover might form in the winter \cite{Abbot08}. If this contribution were to dominate the seasonal variation of the outgoing longwave radiation, then it might shift the minimum of ${\cal F}_l$ from summer toward winter, as in Figure~\ref{fig:Fs-and-Fl}(c).

Mathematically, this study provides an example of how bifurcation behavior of a smoothed system may be reflected in its Filippov counterpart, derived by taking the limit as a smoothing parameter goes to zero, and how the Filippov system may in turn help us understand behavior of the smoothed system. It also gives a unique perspective on the role of repelling sliding intervals in bifurcation behavior of a system with positive feedback that is periodically forced. The link between the width of the sliding intervals and the extent of hysteresis in the system is essential in exploring system behavior. We expect that this link may carry over to other systems with positive feedback, where the feedback strength is periodically modulated. Specifically, we anticipate the perspective used in this case study to apply when the positive feedback involves a sigmoidal function that can be taken to a PWS limit, as we did with the ice-albedo feedback.

\begin{appendices}
\section{Bifurcation set}
\label{app:comparison}

In this appendix we further explore how the bifurcation behavior of the Filippov system (\ref{eq:system-summary}),(\ref{eq:PWSFE12}) compares to that of the original albedo-smoothed system (\ref{eq:dEdt}),(\ref{eq:TE12}),(\ref{eq:FE12}) used in the case study of Section~\ref{sec:case-study}. First we define the bifurcation points for the Filippov system using the existence conditions computed in Section \ref{sec:solutions}. Then we show the bifurcation sets of both systems in the $(\Delta\alpha,L_m)$- and $(\phi,L_m)$-parameter planes, holding all other parameters fixed at the default values, as given in Figure \ref{fig:solution-trajectories}. We compare these results with corresponding figures given in Figure 9(F),(G) in \cite{Eisenman12}. We show that bifurcation behavior of the Filippov system is qualitatively similar to that of the albedo-smoothed system for the parameters studied. 

\begin{figure}[!t]
	\centering
  	\includegraphics[scale=1]{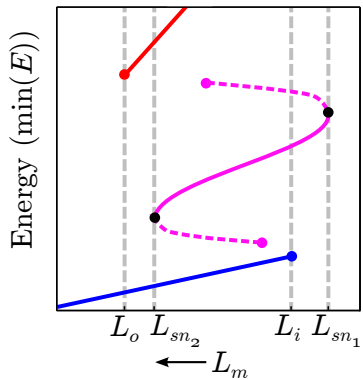}
		\caption{Schematic bifurcation diagram illustrating the locations of the grazing-sliding points $L_{o,i}$ and the saddle-node points $L_{sn_{1,2}}$.}
\label{fig:bif-diag-Lois12}
\end{figure}

For the albedo-smoothed system investigated by Eisenman the only bifurcations identified in \cite{Eisenman12} are saddle-node bifurcations. However, he also keeps track of the existence regions for the three types of periodic states; for the albedo-smoothed problem, these existence boundaries are associated with the $E=0$ boundary: perennially ice-covered (or ice-free) states become seasonally ice-free once their maximum (or minimum) first hits $E=0.$ In the context of nonsmooth systems, solutions which cross a discontinuity boundary as a parameter is changed are no longer ``piecewise-topologically equivalent'' to solutions that remain on one side of the boundary \cite{diBernardo}. The parameter value where this collision with the $E=0$ discontinuity first occurs is a ``grazing'' bifurcation point of the albedo-smoothed system \cite{diBernardo}. For the Filippov system, the existence boundaries are associated with grazing-sliding bifurcations, and there is also the possibility of saddle-node bifurcations on the seasonal solution branches. We label the grazing-sliding bifurcation point of the perennially ice-free (ocean) branch as $L_o$, the grazing-sliding bifurcation point of the perennially ice-covered branch as $L_i,$ and the saddle-node bifurcation with the highest (lowest) minimum energy as $L_{sn_1}$ ($L_{sn_2}$); see Figure \ref{fig:bif-diag-Lois12}.

\begin{figure}[!p] 
	\centering
	\begin{tabular}{ccc}
	\multicolumn{3}{c}{\subfloat[]{\includegraphics[scale=0.96]{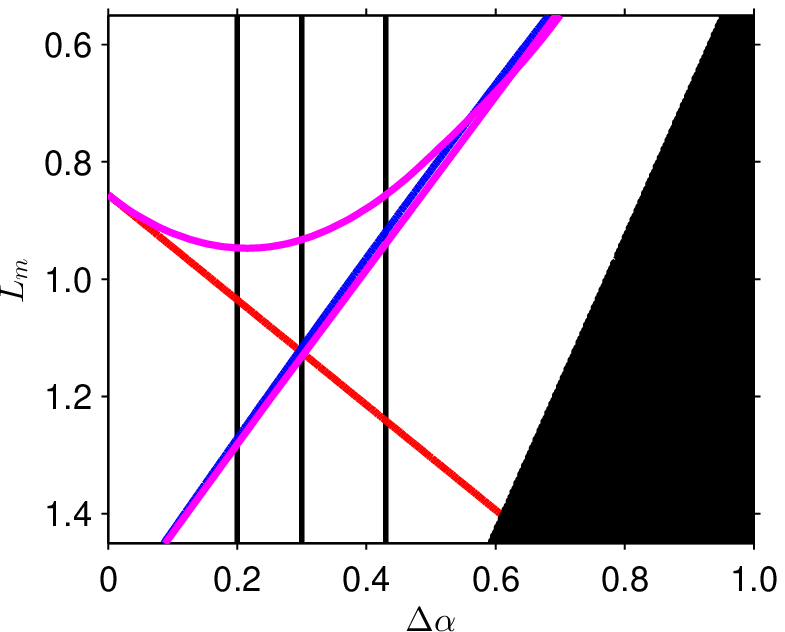}}} \\
	\subfloat[]{\includegraphics[scale=1]{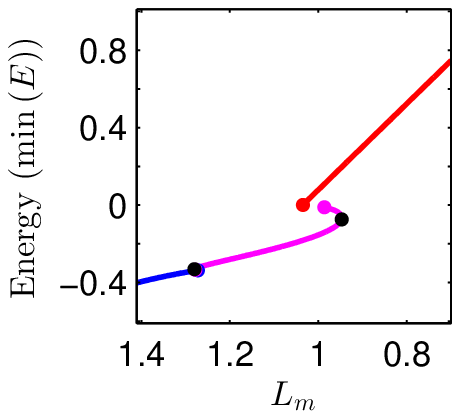}} &
	\subfloat[]{\includegraphics[scale=1]{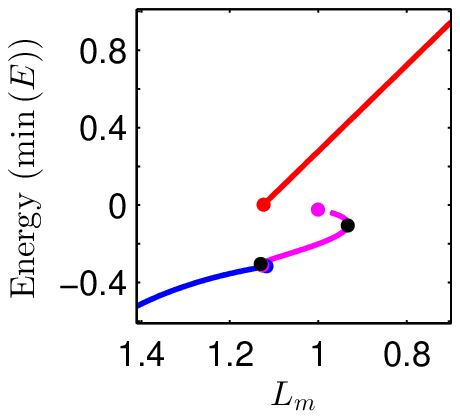}} &
	\subfloat[]{\includegraphics[scale=1]{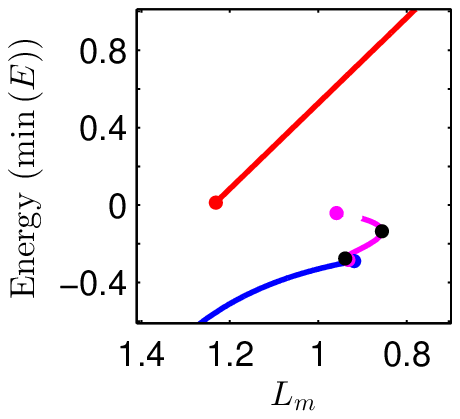}} \\
	\subfloat[]{\includegraphics[scale=1]{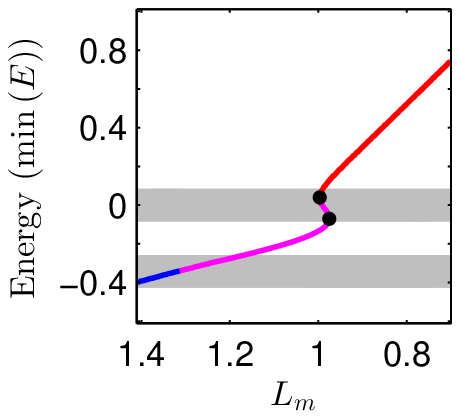}} &
	\subfloat[]{\includegraphics[scale=1]{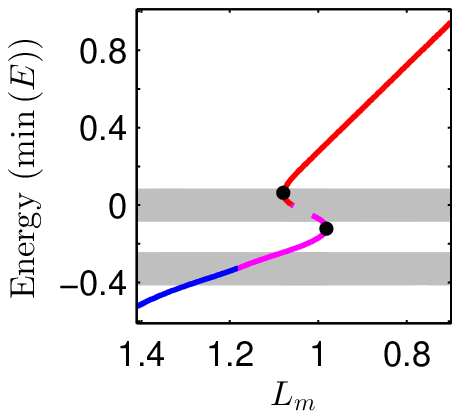}} &
	\subfloat[]{\includegraphics[scale=1]{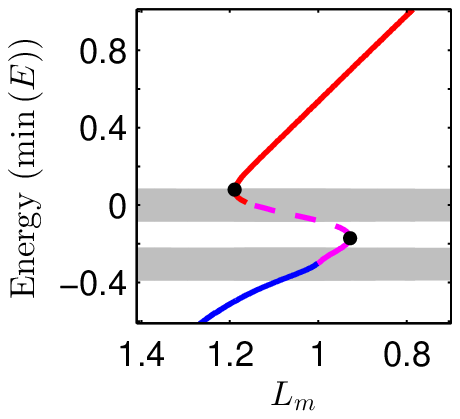}} \\
	\end{tabular}
\caption{(a) Bifurcation set showing the values of $L_{o,i,sn_{1,2}}$ (red, blue, and magenta curves, respectively) of the Filippov system (\ref{eq:system-summary}),(\ref{eq:PWSFE12}) as a function of $\Delta\alpha$, using the default parameters in Figure \ref{fig:solution-trajectories}. Vertical lines at $\Delta\alpha=0.2,$ $0.3,$ and $0.43$ show where bifurcation points occur in the bifurcation diagrams in (b),(c), and (d), respectively. (e),(f), and (g) show bifurcation diagrams of the albedo-smoothed system (\ref{eq:dEdt}),(\ref{eq:TE12}),(\ref{eq:FE12}) for $\Delta E = 0.08$, corresponding to the diagrams for the Filippov system directly above.}
\label{fig:bif-curve-explanation}
\end{figure}

The values of $L_o$ and $L_i$ are determined by the existence boundary for the corresponding perennial solution state, as derived in Section \ref{sec:solutions} ((\ref{eq:ice-free-exist}) and (\ref{eq:ice-cov-exist}), respectively). The saddle-node bifurcation points of the seasonally ice-free branch, $L_{sn_{1,2}}$, simultaneously satisfy (\ref{eq:sn}) and enable us to find $(t_m,t_f)$ which satisfy (\ref{eq:tmtfcondition1}),(\ref{eq:tmtfcondition2}). Note that $L_{o,i}$ always exist, but $L_{sn_{1,2}}$ may collide and annihilate each other in a ``hysteresis point.'' The grazing-sliding bifurcations associated with the existence boundaries for the seasonally ice-free state have no direct counterpart in the albedo-smoothed system, so we omit them for our bifurcation curve study. Figure \ref{fig:bif-curve-explanation}(a) shows the bifurcation set for $L_{o,i,sn_{1,2}}$ as a function of $\Delta\alpha,$ derived from our analysis in Section \ref{sec:solutions}. The black region in (a) corresponds to where there is no repelling sliding interval, and the boundary for this region occurs where $t_b=t_c.$

\begin{figure}[!ht]
	\centering
	\begin{tabular}{ccc}
		\subfloat[]{\includegraphics[scale=1]{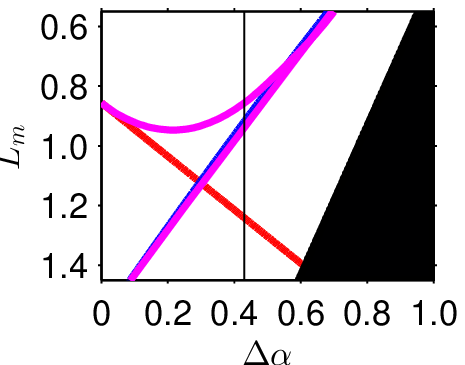}} &
		\subfloat[]{\includegraphics[scale=1]{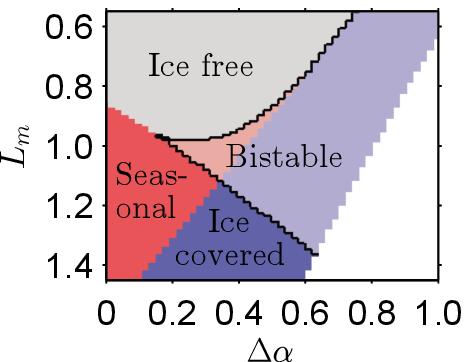}} &
		\subfloat[]{\includegraphics[scale=1]{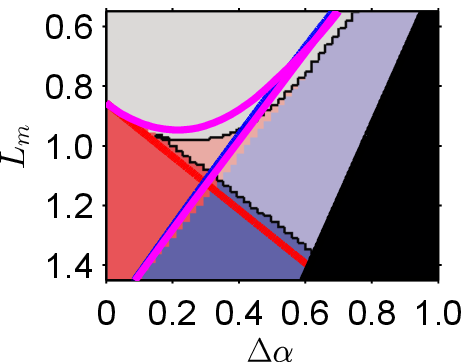}} \\
	\end{tabular}
	\caption{(a) Bifurcation set for the Filippov system (\ref{eq:system-summary}),(\ref{eq:PWSFE12}) as a function of $\Delta\alpha$, using the default parameters. (b) Bifurcation set and classes of stable fixed points of the albedo-smoothed system (\ref{eq:dEdt}),(\ref{eq:TE12}),(\ref{eq:FE12}), reproduced with permission of the author from Figure 9(G) in \cite{Eisenman12}. Black boundaries between different classes represent smooth bifurcations; gray, red, and blue shading is where perennially ice-free, seasonally ice-free, and perennially ice-covered solutions are stable; light red is where both perennially and seasonally ice-free solutions are stable, and light blue is where both perennially ice-covered and ice-free solutions are stable (see \cite{Eisenman12} for further description). (c) Bifurcation set of the Filippov system overlaid on the bifurcation set of the albedo-smoothed system.}
	\label{fig:existence-region}
\end{figure}

Figure \ref{fig:existence-region} presents a comparison of the bifurcation sets for the Filippov and albedo-smoothed systems in the $(\Delta\alpha,L_m)$-parameter plane. In our analysis we use the grazing-sliding points $L_{o,i}$ as proxies for the existence boundaries for the stable perennial branches, as suggested in Figure \ref{fig:bif-diag-comparison}(b) and similar examples in Section \ref{sec:case-study}. Figure \ref{fig:existence-region}(a) shows the bifurcation set for the values of $L_{o,i,sn_{1,2}}$ as a function of $\Delta\alpha$. In (b) we indicate the corresponding bifurcation set of the albedo-smoothed system, which is reproduced from \cite{Eisenman12} and was obtained in~\cite{Eisenman12} by numerically constructing Poincar\'e maps of the albedo-smoothed system. In (c) we overlay the bifurcation set of the Filippov system on the bifurcation set of the albedo-smoothed system for direct comparison. 

We make the following remarks about Figure \ref{fig:existence-region}:
\begin{enumerate}[nolistsep]
\item[1.] The vertical line in (a) indicates the bifurcation points for the default parameters (cf.~the bifurcation diagram of Figure \ref{fig:bif-diag-comparison}(a)).
\item[2.] The bifurcation set for the Filippov system in (a) shows the order and distance between bifurcation points as the additional parameter $\Delta\alpha$ is varied.
\item[3.] The black region in (a) is similar to the white region of the albedo-smoothed system in (b), which for the smoothed system indicates where there is unphysical runaway cooling \cite{Eisenman12}.
\item[4.] In (c), notice that the bifurcation set for the Filippov system is qualitatively similar to that of the albedo-smoothed system. In particular, note that $L_{o,i}$ (the blue and red curves) are very close to the curves that signal the end of the perennial solution branches in the albedo-smoothed system, as illustrated in Figure~\ref{fig:bif-diag-comparison}(b). 
\end{enumerate}

\begin{figure}[!t]
	\centering
	\begin{tabular}{ccc}
		\subfloat[]{\includegraphics[scale=1]{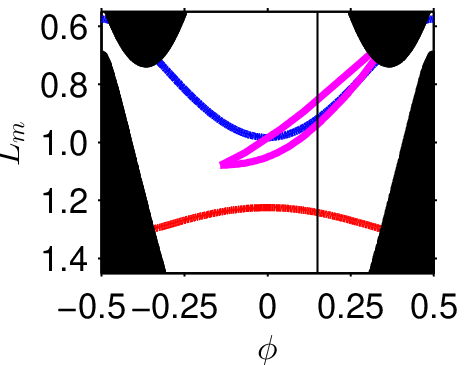}} &
		\subfloat[]{\includegraphics[scale=1]{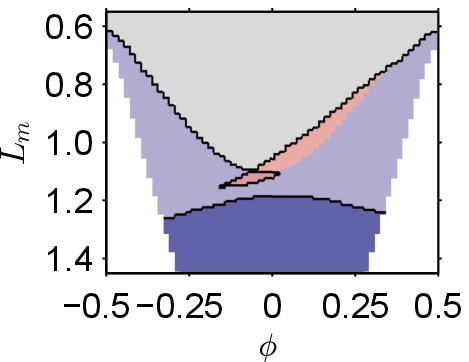}} &
		\subfloat[]{\includegraphics[scale=1]{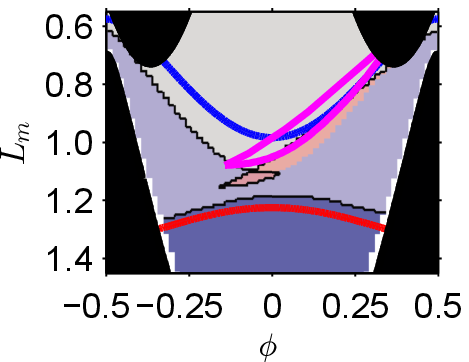}} \\
	\end{tabular}
	\caption{(a) Bifurcation set showing the values of $L_{o,i,sn_{1,2}}$ (red, blue, and magenta curves, respectively) for the Filippov system as a function of $\phi$, using the default parameters. (b) Bifurcation set and classes of stable fixed points of the albedo-smoothed system, as in Figure 9(F) in \cite{Eisenman12}; shading is as described in Figure \ref{fig:existence-region}. (c) Bifurcation set of the Filippov system overlaid on the bifurcation set of the albedo-smoothed system.}
	\label{fig:phi-existence-region}
\end{figure}

In Figure \ref{fig:phi-existence-region}, we provide a similar comparison of the bifurcation set between the Filippov and albedo-smoothed systems, for varying $\phi.$ Note that our analysis requires us to exclude a region of parameter space not excluded in the albedo-smoothed system, due to the presence of attracting sliding intervals in the Filippov system. Also note that in (b) we have shifted the range of the phase $\phi$ from the original figure in \cite{Eisenman12} so that $\phi\in(-0.5,0.5].$ In general, while there are noticeable quantitative differences, there is significant qualitative similarity, as in Figure \ref{fig:existence-region}.

\section{Inverse mapping}
\label{app:inverse-mapping}

In this appendix we determine an inverse mapping, when it exists, between Eisenman's  physical parameters, $\Delta\alpha$, $S_a$, $L_m$, $L_a$, $\phi$ of (\ref{eq:PWSFE12}), and the parameters we introduced for our analysis in Section~\ref{sec:simplified-model}, namely, $\overline{F_\pm}$, $\widetilde{F_\pm}$ and $\Delta\psi\equiv\psi_+-\psi_-$ of~(\ref{eq:simple-FE12}). We are interested in the mapping so that, for example, we can construct  a path through Eisenman's five-parameter space  ($L_m,$ $\Delta\alpha,$ $S_a,$ $L_a,$ $\phi$) that is associated with varying the relative phase $\Delta\psi,$ while holding the four amplitudes ($\overline{F_+},$ $\overline{F_-},$ $\widetilde{F_+},$ $\widetilde{F_-}$) fixed. At the end of this appendix is a table of parameter sets used to generate the bifurcation diagrams in Figure~\ref{fig:vary-DPsi-and-Fpmtilde}.

Here we reproduce the two expressions that are relevant to our construction:
\begin{equation}\label{eq:PWSFE12-appendix}
	F_\pm(\tau) = (1 \pm \Delta\alpha)(1-S_a \cos(2 \pi \tau)) - (L_m + L_a \cos(2 \pi (\tau-\phi)))
\end{equation}
and 
\begin{equation}\label{eq:trig}
	F_\pm(\tau)=\overline{F_\pm}+\widetilde{F_\pm}\cos(2\pi \tau-\psi_\pm).
\end{equation}
We can relate the phases ($\psi_\pm$) and amplitudes ($\overline{F_\pm}, \widetilde{F_\pm}$) to Eisenman's physical parameters via
\begin{eqnarray}\label{eq:amplitudes}
	\overline{F_\pm} &=& 1-L_m\pm\Delta\alpha, \nonumber\\
	\widetilde{F_\pm} \cos(\psi_\pm) &=& -(1\pm\Delta\alpha)S_a-L_a\cos(2\pi\phi), \\
	\widetilde{F_\pm}\sin(\psi_\pm) &=& -L_a\sin(2\pi\phi).\nonumber
\end{eqnarray}
We introduce the following shorthand for some of the relationships in (\ref{eq:amplitudes}):
\begin{eqnarray}\label{eq:newones}
	\widetilde{F_+}\cos(\psi_+)&=&-S_+-L_{ac}, \nonumber\\
	\widetilde{F_-}\cos(\psi_-)&=&-rS_+-L_{ac}, \\
	\widetilde{F_\pm}\sin(\psi_\pm)&=&-L_{as},\nonumber
\end{eqnarray}
where
\begin{align}\label{eq:defs}
	\begin{aligned}
		S_+&\equiv (1+\Delta\alpha)S_a\in[S_a, 2S_a], \\
		r&\equiv \frac{1-\Delta\alpha}{1+\Delta\alpha}\in [0,1], \\
		L_{ac}&\equiv L_a\cos(2\pi\phi)\in[-L_a,L_a], \\
		L_{as}&\equiv L_a\sin(2\pi\phi),\ {\rm and \ where}\ L_{ac}^2+L_{as}^2=L_a^2.
	\end{aligned}
\end{align}

Here we enumerate some preliminary steps used in constructing an inverse mapping:
\begin{enumerate}[nolistsep]
\item{}  $\overline{F_+}$ and $\overline{F_-}$  determine both $L_m$ and $\Delta\alpha$, which in turn determine $r$ in (\ref{eq:newones})--(\ref{eq:defs}):
\begin{align}\label{eq:overlineFpm}
	\begin{aligned}
		\Delta\alpha&=\frac{1}{2}\Bigl(\overline{F_+}-\overline{F_-}\Bigr), \\
		L_m&=1-\frac{1}{2}\Bigl(\overline{F_+}+\overline{F_-}\Bigr).
	\end{aligned}
\end{align}
Note that the conditions $\Delta\alpha\in[0,1]$ and $L_m\ge 0$ place some restrictions on the possible values of $\overline{F_\pm}$.

\item{} We next use the values of $\widetilde{F_+}$ and $\widetilde{F_-}$ from (\ref{eq:newones}) to determine the parameter combinations $S_+$ and $L_{ac}$, where
\begin{eqnarray}
	\widetilde{F_+}^2&=&L_a^2+S_+^2+2S_+L_{ac}, \nonumber\\
	\widetilde{F_-}^2&=&L_a^2+r^2S_+^2+2rS_+L_{ac}. \nonumber
\end{eqnarray}
This pair of equations can be solved, given $\widetilde{F_\pm}$ and $r$ determined in step 1, to determine $S_+$ and $L_{ac}$ as functions of $L_a^2$. Specifically, we find
\begin{align}\label{eq:params}
	\begin{aligned}
		S_+^2&=\frac{r\widetilde{F_+}^2-\widetilde{F_-}^2+L_a^2(1-r)}{r(1-r)}, \\
		S_+ L_{ac}&=
		\frac{r^2\widetilde{F_+}^2-\widetilde{F_-}^2+L_a^2(1-r^2)}{2r(r-1)}.
	\end{aligned}
\end{align}
The positive root of the $S_+^2$ equation is used to determine  $S_+\ge 0$. Note that  there are some restrictions placed on the possible values for $\widetilde{F_\pm}$ by the physically allowed ranges on $S_+^2$ and $L_{ac}$.

\item{} We next obtain an equation for $L_a$ in terms of the quantities calculated in the previous steps and $\tan(\Delta\psi)$, where
	\begin{align}\label{eq:tandpsi}
	\tan(\Delta\psi)&=\frac{\sin(\psi_+)\cos(\psi_-)-\cos(\psi_+)\sin(\psi_-)}{\cos(\psi_+)\cos(\psi_-)+\sin(\psi_+)\sin(\psi_-)}\nonumber\\
	&=\frac{(r-1)S_+L_{as}}{L_a^2+rS_+^2+S_+L_{ac}(1+r)}.
	\end{align}
The values of $S_+^2$ and $S_+L_{ac}$ are given by~(\ref{eq:params}) in terms of the amplitudes $\widetilde{F_\pm}$ and $L_a^2$. Also, $L_{as}=\pm\sqrt{L_{a}^2-L_{ac}^2}$, where the choice of the sign will ultimately be determined by the sign of $\tan(\Delta\psi)$, and in turn that will determine the sign of $\phi$ since $L_{as}=L_a\sin(2\pi\phi)$, where $\phi\in(-\frac{1}{2},\frac{1}{2}]$.  Combining all of this, we find that $L_a^2$ solves the following quadratic in $\widetilde{F_+},\,\widetilde{F_-},\,r$, and $\tan^2(\Delta\psi)$:
\begin{equation}\label{eq:La2quad}
	L_a^4-\frac{2(\widetilde{F_-}^2+r^2\widetilde{F_+}^2)}{(1-r)^2}L_a^2+\frac{(\widetilde{F_-}^2-r^2\widetilde{F_+}^2)^2+(\widetilde{F_-}^2+r^2\widetilde{F_+}^2)^2\tan^2(\Delta\psi)}{(1-r)^4(1+\tan^2(\Delta\psi))}=0.
\end{equation}
Again, there are restrictions on the values of $\widetilde{F_\pm}$ and $\Delta\psi$ in order for this to have real solutions. Moreover, there can be two positive, real solutions $L_a^2$ to this equation. We find, when there are two distinct solutions for $L_a^2,$ that one is much larger than the other, and it is typically the smaller one that is of interest to us. 

\end{enumerate}

To summarize, given a set of parameters $\overline{F_\pm}$, $\widetilde{F_\pm}$ and $\Delta\psi$ of~(\ref{eq:trig}), we compute parameters for the Eisenman model~(\ref{eq:PWSFE12-appendix}) as follows:
\begin{enumerate}[nolistsep]
\item{} Compute $\Delta\alpha$ and $L_m$ using~(\ref{eq:overlineFpm}), and then compute $r$ using its definition in~(\ref{eq:defs}).
\item{} Compute $L_a^2$ from the quadratic~(\ref{eq:La2quad}), which may have two roots. The choice of root determines $L_a$, which is nonnegative.
\item{} Compute $S_+$ and $L_{ac}$ from~(\ref{eq:params}). $S_+$ and $\Delta\alpha$ determine $S_a$ using~(\ref{eq:defs}). 
\item{}  $L_{as}$ can be determined, up to a sign, from $L_a^2$ and $L_{ac}$. Then the sign of $L_{as}$ is determined from the sign of $\tan(\Delta\psi)$ in~(\ref{eq:tandpsi}). Finally the phase $\phi$ is determined from $L_{ac}$ and $L_{as}$, defined in~(\ref{eq:defs}).
\end{enumerate}

\begin{table}[!ht]
\centering
\setlength{\tabcolsep}{1em} 
{\renewcommand{\arraystretch}{1.5}
\begin{tabular}{|l|c|c|c||c|c|c|}
	\hline 
	 & \multicolumn{3}{c||}{\textbf{Mapped parameters}} & \multicolumn{3}{c|}{\textbf{Original parameters}}\tabularnewline
	\hline 
	$\vphantom{\widetilde{F_{+}}}$ & $\Delta\psi$ & $\widetilde{F_{+}}$ & $\widetilde{F_{-}}$ & $S_{a}$ & $L_{a}$ & $\phi$\tabularnewline
	\hline 
	\multirow{4}{*}{\textbf{Varying $\Delta\psi$}} & $-0.21$ & $2.64$ & $1.41$ & $1.50$ & $0.73$ & $\hphantom{-}0.15$\tabularnewline
	\cline{2-7} 
	 & $\hphantom{-}0.00$ & $2.64$ & $1.41$ & $1.43$ & $0.60$ & $\hphantom{-}0.00$\tabularnewline
	\cline{2-7} 
	 & $\hphantom{-}0.30$ & $2.64$ & $1.41$ & $1.58$ & $0.85$ & $-0.20$\tabularnewline
	\cline{2-7} 
	 & $\hphantom{-}0.51$ & $2.64$ & $1.41$ & $1.82$ & $1.19$ & $-0.28$\tabularnewline
	\hline 
	\hline 
	\multirow{4}{*}{\textbf{Varying $\widetilde{F_{+}}$}} & $-0.21$ & $1.45$ & $1.41$ & $0.34$ & $1.42$ & $\hphantom{-}0.27$\tabularnewline
	\cline{2-7} 
	 & $-0.21$ & $2.00$ & $1.41$ & $0.79$ & $1.09$ & $\hphantom{-}0.14$\tabularnewline
	\cline{2-7} 
	 & $-0.21$ & $2.64$ & $1.41$ & $1.50$ & $0.73$ & $\hphantom{-}0.15$\tabularnewline
	\cline{2-7} 
	 & $-0.21$ & $4.00$ & $1.41$ & $3.06$ & $0.59$ & $\hphantom{-}0.38$\tabularnewline
	\hline 
	\hline 
	\multirow{4}{*}{\textbf{Varying $\widetilde{F_{-}}$}} & $-0.21$ & $2.64$ & $0.80$ & $2.17$ & $0.52$ & $\hphantom{-}0.43$\tabularnewline
	\cline{2-7} 
	 & $-0.21$ & $2.64$ & $1.41$ & $1.50$ & $0.73$ & $\hphantom{-}0.15$\tabularnewline
	\cline{2-7} 
	 & $-0.21$ & $2.64$ & $2.00$ & $0.93$ & $1.65$ & $\hphantom{-}0.15$\tabularnewline
	\cline{2-7} 
	 & $-0.21$ & $2.64$ & $2.50$ & $0.63$ & $2.47$ & $\hphantom{-}0.25$\tabularnewline
	\hline 
\end{tabular}
}
\caption{Approximate parameters from~(\ref{eq:PWSFE12}) used to generate the bifurcation diagrams in Figure~\ref{fig:vary-DPsi-and-Fpmtilde}, calculated using the inverse mapping described in this appendix. The remaining parameters of the full model ($B$, $\zeta,$ and $\Delta\alpha$) have been set to the values given in Figure~\ref{fig:solution-trajectories}.}
\label{tab:bif-pars}
\end{table}

\end{appendices}

\section*{Acknowledgments}
The authors would like to acknowledge support of the Mathematics and Climate Research Network (MCRN) under NSF grants DMS-0940261 and DMS-0940262. KH was additionally supported by a Microsoft Research Graduate Women's Scholarship and an NSF Graduate Research Fellowship under Grant No. DGE-1324585.


\end{document}